\documentclass[11pt,titlepage]{article}
\usepackage{latexsym}
\usepackage{amsfonts}
\usepackage{amsmath}
\usepackage{amssymb}
\usepackage{graphicx}
\usepackage{tikz}
\usetikzlibrary{arrows, decorations.markings}
\newcommand\blackslug{\hbox{\hskip 1pt \vrule width 4pt height 8pt depth 1.5pt
        \hskip 1pt}}
\newcommand\bbox{\hfill \quad \blackslug \medbreak}

\newtheorem{theorem}{}[section]

\newtheorem{conjecture1}{Conjecture}[section]
\newtheorem{lemma}{Lemma}
\newtheorem{definition}{Definition}
\newcommand{\Proof}{\noindent{\bf Proof.}\ \ }

\tikzstyle{every node}=[circle, draw, fill=black!50, inner sep=0pt, minimum width=20pt]

\tikzstyle{input}=[circle,
                                    thick,
                                    minimum size=2.2cm,
                                    draw=black!80,
                                    fill=white!20]

\tikzstyle{matrx}=[rectangle,
                                    thick,
                                    minimum size=1.8cm,
                                    draw=black!80,
                                    fill=white!20]

\tikzstyle{vecArrow} = [thick, decoration={markings,mark=at position
   1 with {\arrow[semithick]{open triangle 60}}},
   double distance=1.8pt, shorten >= 5.5pt,
   preaction = {decorate},
   postaction = {draw,line width=1.4pt, white,shorten >= 4.5pt}]
\tikzstyle{innerWhite} = [semithick, white,line width=1.4pt, shorten >= 4.5pt]

\title{All known prime Erd\H{o}s-Hajnal tournaments satisfy $\epsilon(H) = \Omega(\frac{1}{|H|^{5}\log(|H|)})$}
\author{Krzysztof Choromanski,
Columbia University\\
New York, USA
}
\date{August 1, 2014; revised \today}

\newcommand{\Keywords}{the Erd\H{o}s-Hajnal Conjecture, transitive subtournaments, Ramsey theory, prime tournaments, galaxies}

\begin{document}
\maketitle
\begin{abstract}
We prove that there exists $C>0$ such that $\epsilon(H) \geq \frac{C}{|H|^{5}\log(|H|)}$, where $\epsilon(H)$ is the Erd\H{o}s-Hajnal coefficient of the tournament $H$, for every prime tournament $H$ for which the celebrated Erd\H{o}s-Hajnal Conjecture has been proven so far. This is the first polynomial bound on the EH coefficient obtained for all known prime Erd\H{o}s-Hajnal tournaments, in particular for infinitely many prime tournaments.
As a byproduct of our analysis, we answer affirmatively the question whether there exists an infinite family of prime tournaments $H$ with $\epsilon(H)$ lower-bounded by $\frac{1}{\textit{poly}(|H|)}$, where $\textit{poly}$ is a polynomial function.
Furthermore, we give much tighter bounds than those known so far for the EH coefficients of tournaments without large homogeneous sets.
This enables us to significantly reduce the gap between best known lower and upper bounds for the EH coefficients of tournaments.
As a corollary we prove that every known prime Erd\H{o}s-Hajnal tournament $H$ satisfies:
$-5 + o(1) \leq \frac{\log(\epsilon(H))}{\log(|H|)} \leq -1 + o(1)$. No lower bound on that expression was known before.
We also show the applications of those results to the tournament coloring problem. In particular, we prove that for every known prime Erd\H{o}s-Hajnal tournament $H$ every $H$-free tournament has \textit{chromatic number} at most
$O(n^{1-\frac{C}{|H|^{5}\log(|H|)}}\log(n))$, where $C>0$ is some universal constant.
The related coloring can be constructed algorithmically in the quasipolynomial time by following straightforwadly the proof of our main result.
In comparison, the standard Ramsey theory gives only $O(\frac{n}{\log(n)})$ bounds for the tournament chromatic number.
\end{abstract}

\maketitle {\bf Keywords:} \Keywords

\newpage

\section{Introduction}

We focus in this paper on estimating one of the most interesting graph invariants in modern Ramsey graph theory, 
the so-called \textit{EH coefficient} of the tournament (also known as the \textit{Erd\H{o}s-Hajnal coefficient}). The EH coefficient comes from one of the most challenging and still open problems in Ramsey graph theory - the \textit{Erd\H{o}s-Hajnal Conjecture}. 
The Conjecture states that for every tournament $H$ there exists $\epsilon(H) >0$ (the EH coefficient) such that every $n$-vertex $H$-free tournament contains a transitive subtournament of order at least $n^{\epsilon(H)}$. 
Despite many attempts, the Conjecture has been unsolved for more than twenty years now and was proven so far only for some specific classes of tournaments $H$. However derived lower bounds on the EH coefficients $\epsilon(H)$ 
were extremely small for most of them - at most inversely proportional to the Szemeredi tower function. 
The author of this paper is not aware of any published result that proposed bounds not relying on the Szemeredi Lemma. 
On the other hand, the best known upper bounds on the EH coefficient, obtained by the probabilistic method, are of order $O(\frac{1}{|H|^{1-o(1)}})$. 
Thus, like in the case of the \textit{Ramsey number}, the best known lower bounds on the EH coefficient were far away from the best known upper bounds for most of the tournaments. 
We show that for all prime tournaments for which lower bounds greater 
than $0$ were obtained those bounds can be in fact strengthened to be of order 
$\Omega(\frac{1}{|H|^{5}\log(|H|)})$. 
Thus we get polynomial bounds and surprisingly, with very small polynomial degree that
does not depend on the structure of the tournament. At the same time we significantly 
reduce the gap between best known lower and upper bounds on the EH coefficients 
of tournaments. This is the first polynomial bound that works for all of those tournaments and a step towards answering Erd\H{o}s question how the EH coefficients depend on the order of the forbidden tournament $H$.  Prime tournaments play crucial role in the study of the Conjecture and in many other graph theory problems - if the Conjecture is true for prime tournaments then it is true for all the tournaments.
As a corollary we prove new lower bounds on EH coefficients for all tournaments for which the Conjecture was proven so far
and answer affirmatively another question regarding the conjecture: whether there exists a universal
polynomial lower bound for the EH coefficient for some infinite family of prime tournaments.
Our result tightens lower and upper bounds for EH coefficients of several classes of tournaments. For some of them we obtain even tighter bounds. For the family of prime \textit{stars} we are able to prove that $\epsilon(H) = \frac{1}{|H|^{1+o(1)}}$.  
This is the first result giving asympotically tight lower and upper bounds on the EH coefficient for an infinite
family of prime tournaments.
We also give tighter lower bounds on EH coefficients for tournaments $H$ without large homogeneous sets. 

Our results lead to new purely combinatorial coloring algorithms for classes of graphs characterized by forbidden patterns.
These classes play an important role in graph theory. For instance, every graph with the topological ordering of vertices can be 
equivalently described as not having directed cycles and every transitive tournament - as not having directed triangles.
A finite graph is planar if and only if it does not contain $K_{5}$ (the complete graph on five vertices) or $K_{3,3}$ (complete bipartite graph on six vertices with
two equal-length color classes) as a minor. One of the deepest results in graph theory, the Robertson-Seymour theorem (\cite{robertson}), states that every family of graphs 
(not necessarily planar graphs) that is closed under minors can be defined by a finite set of forbidden minors. These classes include: forests, pseudoforests, linear forests
(disjoint unions of path graphs), planar and outerplanar graphs, apex graphs, toroidal graphs, graphs that can be embedded on the two-dimensional manifold, graphs with
bounded treewidth, pathwidth or branchwidth and many more. This theorem has also a directed version. 
The other examples include classes of graphs that can be colored with significantly fewer than $\Omega(\frac{n}{\log(n)})$ colors (for instance the classes of graphs that are
$c$-colorable for some constant $c>0$ that were intensively studied before). All those classes can be described as not having some nontrivial forbidden structures (either induced
subgraphs in the undirected scenario or subtournaments in the directed setting).
Thus classes of graphs described by forbidden patterns appear very often
in both directed and undirected setting. We should notice that not having a certain
graph as a minor is a much more restrictive assumption than not having a certain graph $H$ as an induced subgraph. That is why our setting is much more general.
We present our methods in the directed scenario but they can be translated into undirected one. In fact, as we will see soon, the Conjecture has an equivalent undirected version,
where forbidden patterns $H$ are undirected graphs. Surprisingly, both versions are equivalent.
From what we have said before, it is clear that coloring $H$-free graphs is an important algorithmic problem. Our result is a first step to obtain purely combinatorial
nontrivial coloring algorithms for these classes of graphs.
It is worth to mention here that graph coloring problem is NP-hard even to approximate within mutliplicative factor $n^{1-\epsilon}$ for an arbitrary
fixed $\epsilon>0$. Thus algorithms achieving this for some special classes of graphs are of great interest. 
In our setting $\epsilon$ corresponds to the EH coefficient and that
establishes an intriguing connection between the Conjecture (namely, algorithmic proofs of the lower bounds on EH coefficients) and hardness results for approximating algorithms designed to color graphs.
Among some of the most important known results regarding coloring graphs are algorithms for coloring $3$-colorable graphs with at most $n^{\delta}$ colors for 
some $0<\delta<1$ (\cite{blum}, \cite{karger}, \cite{arora}). 
Most of them rely on noncombinatorial approach such as SDP.
Notice that if $H$ is a graph with no stable sets of order at least $\frac{|H|}{3}$ (such a graph can be easily constructed
randomly) then every $3$-colorable graph is $H$-free (and as mentioned earlier, this observation can be easily generalized to the $c$-colorable graphs).
Therefore all those considered and heavily investigated classes of graphs are captured by the forbidden pattern framework.

\subsection{Notation}
We use $||$ to denote the size of the set. Let $G$ be a graph. We denote by $V(G)$ the set 
of its vertices and by $E(G)$ the set of its edges. Sometimes instead of writing $|V(G)|$ we use 
shorter notation $|G|$. We call $|G|$ the \textit{size of $G$} (or \textit{order of $G$}). 
For a subset $S \subseteq V(G)$ we denote by $G|S$ the subgraph of $G$ induced by $S$. 
A \textit{clique} in an undirected graph is a set of pairwise adjacent vertices. 
An \textit{independent set} in the undirected graph is a set of pairwise nonadjacent vertices.

A \textit{tournament} is a directed graph  such that for every pair $v$ and $w$ of vertices, exactly one of the edges $(v,w)$ or $(w,v)$ exists. If $(v,w)$ is an edge of the tournament then we say that $v$ is $\textit{adjacent to}$ $w$ and $w$ is \textit{adjacent from} $v$. 
A subset $S_{1} \subseteq V(T)$ of the vertices of a tournament $T$ is \textit{adjacent to} another subset 
of the vertices $S_{2} \subseteq V(T)$ if every vertex of $S_{1}$ is adjacent to every vertex of $S_{2}$.
Then we also say that $S_{2}$ is adjacent from $S_{1}$.
The \textit{indegree} of a vertex $v$ of a tournament $T$ is the number of vertices $w \in V(T)$ such that $(w,v) \in E(T)$.
Similarly, the \textit{outdegree} of a vertex $v$ of a tournament $T$ is the number of vertices $w \in V(T)$ such that $(v,w) \in E(T)$. 
A tournament is \textit{transitive} if it contains no directed cycle. For the set of vertices $V=\{v_{1},v_{2},...,v_{k}\}$ we say that an ordering $(v_{1},v_{2},...,v_{k})$ is \textit{transitive} if $v_{1}$ is adjacent to all other vertices of $V$, $v_{2}$ is adjacent to all other vertices of $V$ but $v_{1}$, etc. A subset $S \subseteq V(T)$ is \textit{transitive} if it induces a transitive tournament. 
For a tournament $H$ we say that a tournament $T$ is $H$-free if $T$ does not contain $H$ as an induced subtournament.
We denote by $C_{5}$ a unique tournament on five vertices where every vertex has indegree $2$. Other tournaments such as stars and galaxies will be defined later in the paper.
We denote by $tr(T)$ the size of the largest transitive subtournament of a tournament $T$.
A \textit{coloring} of the tournament $T$ is a partitioning of the set of its vertices into transitive subsets. 
A partitioning with minimal number of parts is called a \textit{chromatic number} $\chi(T)$ of the tournament $T$.

All logarithms used in the paper are natural.

\subsection{The Conjecture and the \textit{substitution procedure}}

A celebrated unresolved Conjecture~ of Erd\H{o}s and Hajnal states that:  

\begin{conjecture1}
\label{EHCon}
For every tournament $H$ there exists $\epsilon(H)>0$ such that every $n$-vertex $H$-free tournament contains a transitive subtournament of size at least $n^{\epsilon(H)}$.
\end{conjecture1}

In fact the Conjecture was first proposed in the undirected version by Erd\H{o}s and Hajnal but was proven to have an 
equivalent directed version above by Alon, Pach and Solymosi in 2001 (see: \cite{alon}). 
The undirected version (see: \cite{erdos0}) states that:

\begin{conjecture1}
\label{EHCon-undirected}
For every undirected graph $H$ there exists $\epsilon(H)>0$ such that every $n$-vertex graph $G$ that does not contain $H$ as an induced subgraph contains a clique or an independent set of size at least $n^{\epsilon(H)}$.
\end{conjecture1}

If for a given tournament $H$ there exists $\epsilon(H)>0$ then we say that \textit{$H$ satisfies the Erd\H{o}s-Hajnal Conjecture with $\epsilon(H)$} or simply: \textit{$H$ satisfies the Erd\H{o}s-Hajnal Conjecture}. 
Sometimes we say that H is the \textit{ Erd\H{o}s-Hajnal tournament} or simply: EH tournament.
The coefficient $\epsilon(H)$ in the statement is called the \textit{EH coefficient}
(or \textit{the Erd\H{o}s-Hajnal coefficient}).
\\

For any tournament $H$ with vertex set $V(H)=\{v_{1},...,v_{h}\}$ and
tournaments $F_{1},...,F_{h}$ let $H(F_{1},...,F_{h})$ denote the tournament obtained from $H$ by replacing each $v_{i}$ with a copy of $F_{i}$, and making
a vertex of the copy of $F_{i}$ outadjacent to a vertex of a copy of $F_{j}$, $j \neq i$, if and only if $(v_{i},v_{j}) \in E(H)$. The copies of
$F_{i}$, $i=1,...,h$, are assumed to be vertex disjoint. The procedure of constructing $H(F_{1},...,F_{h})$ from $H, F_{1},...,F_{h}$ is called the \textit{substitution procedure}.

A subset of vertices $S \subseteq V(H)$ of a tournament $H$ is called \textit{homogeneous} if for every $v \in V(H) \backslash S$ the following holds: either $\forall_{w \in S} (w,v) \in E(H)$ or $\forall_{w \in S} (v,w) \in E(H)$. A homogeneous set $S$ is called \textit{nontrivial} if $|S|>1$ and $S \neq V(H)$. A tournament is called \textit{prime} if it does not have nontrivial homogeneous sets. Alon, Pach and Solymosi proved that if the Conjecture is false, then the smallest counterexample is prime. 
They did it by showing an upper bound on $\epsilon(H(F_{1},..., F_{h}))$ as the function of $\epsilon(F_{1}),..., \epsilon(F_{h})$ (see: \cite{alon}). 
A \textit{homogeneous partitioning} of the set of vertices of a tournament is a partitioning of its 
vertices into homogeneous sets. A homogeneous partitioning is called \textit{nontrivial} if it does
not consists of just one set (the set of all the vertices of the tournament). The partitioning number
$p(H)$ of a tournament $H$ is the smallest possible number of parts in the nontrivial homogeneous partitioning. 
Note that if $\{V_{1},V_{1},...,V_{r}\}$ is a homogeneous partitioning of $V(H)$ for $1 \leq i < j \leq r$ then
either $V_{i}$ is adjacent to $V_{j}$ or $V_{j}$ is adjacent to $V_{i}$.

\section{Main results and related work}

We show in this paper that there exists $C>0$ such that $\epsilon(H) \geq \frac{C}{|H|^{5}\log(|H|)}$ for every prime tournament $H$ for which the Conjecture is known.
More precisely, we show the following:

\begin{theorem}
\label{conjecturetheorem}
There exists $C>0$ such that if $H$ is a prime galaxy then 
$$\epsilon(H) \geq\frac{C}{|H|^{5}\log(|H|)}.$$
\end{theorem}

This is the first polynomial bound on the EH coefficients for all prime tournaments for which the Conjecture has been proven so far.
It is enough to focus on galaxies since the only prime tournaments for which the Conjecture
was proven to be true so far are prime galaxies, tournament $C_{5}$ and two more six-vertex tournaments.
In fact a randomly chosen galaxy is prime with high probability.

As an immediate corollary of Theorem \ref{conjecturetheorem} and the very well-known fact that every prime
tournament $H$ satisfies: $\epsilon(H) \leq C\frac{\log(|H|)}{|H|}$ for some universal constant $C$ (see: Appendix A) , we get:

\begin{theorem}
 \label{tighttheorem}
 Every prime galaxy $H$ satisfies: 
 $-5 + o(1) \leq \frac{\log(\epsilon(H))}{\log(|H|)} \leq -1 + o(1)$.
\end{theorem}

That, according to our previous remarks, gives us:  $-5 + o(1) \leq \frac{\log(\epsilon(H))}{\log(|H|)} \leq -1 + o(1)$ for every known prime  Erd\H{o}s-Hajnal tournament.

Another corollary of our core result is the algorithm for coloring $H$-free tournaments. 
For a fixed Erd\H{o}s-Hajnal tournament $H$ we show that these 
are $O(n^{1-\frac{C}{|H|^{5}\log(|H|)}}\log(n))$-colorable for some universal constant $C>0$.

We are also the first to prove tight asymptotic lower and upper bounds on Erd\H{o}s-Hajnal coefficients
for some infinite classes of prime tournaments. In particular, we prove that:

\begin{theorem}
\label{startheorem}
If $H$ is a prime star then $\epsilon(H) = \frac{1}{|H|^{1+o(1)}}$.
\end{theorem}

As a corollary of our techniques, we prove the following:

\begin{theorem}
 \label{generaltheorem}
 All known Erd\H{o}s-Hajnal tournaments $H$ satisfy: $$\epsilon(H) \geq \frac{C}{e^{(|H|+1)\log(|H|+1)}},$$
 for some constant $C > 0$.
\end{theorem}

The bounds given above for all known Erd\H{o}s-Hajnal tournaments are substantially better than the best
previously known. As a corollary, we  answer affirmatively the question whether there exists an infinite family of prime tournaments $H$ with $\epsilon(H)$ 
lower-bounded by $\frac{1}{\textit{poly}(|H|)}$, where $\textit{poly}$ is a polynomial function.\\

The following theorem turns out to be extremely useful in obtaining strong lower bounds on EH coefficients for
prime tournaments $H$, not only in the context of this paper. It is also interesting in itself.

\begin{theorem}
 \label{largeenoughtheorem}
 Assume that $H$ is a prime tournament and that every $n$-vertex $H$-free tournament contains a 
 transitive subtournament of order at least $c(H) n^{\epsilon(H)}$ for some $c(H), \epsilon(H) > 0$.
 Then every $n$-vertex $H$-free tournament contains a transitive subtournament of order at least $n^{\epsilon(H)}$.
\end{theorem}

The theorem says that, surprisingly, in order to prove that a certain expression is the lower bound on the
EH coefficient of a prime tournament it suffices to prove a similar result where an additional constant $c(H)>0$
is being added as a multiplicative factor to the expression on the size of the transitive subtournament.
This, as we will see very soon, simplifies the analysis very much.

Theorem \ref{generaltheorem} can be used to obtained stronger than the best known so far upper bounds on the chromatic numbers of tournaments defined by forbidden patterns.
We will prove the following.

\begin{theorem}
\label{algorithmtheorem}
For any prime galaxy $H$ there exists a quasipolynomial algorithm that finds a coloring of the $n$-vertex $H$-free 
tournament with $O(n^{1-\frac{C}{|H|^{5}\log(|H|)}}\log(n))$ colors, 
where $C$ is some universal constant.
\end{theorem}

We also significantly improve upper bounds on EH coefficients for tournaments that even though, not necessarily
prime, have relatively small homogeneous sets. By doing it we reduce the gap between lower and upper bounds on EH coefficients
for many more classes of tournaments. We prove that:

\begin{theorem}
 \label{uppertheorem}
  For every $\eta >0$ there exists $C(\eta)>0$ such that every tournament $H$ with the largest nontrivial 
  homogeneous sets of size at most $\frac{\sqrt{h}}{2}$ satisfies:
  $$\epsilon(H) \leq C(\eta) \frac{\log(p(H))}{p(H)^{\frac{1}{2}-\eta}},$$
  where $p(H)$ is the partitioning number of $H$ and $h=|H|$.
\end{theorem}

This result is a significant improvement since the best previously known bounds for the tournaments that
are not necessarily prime were only inversely proportional to the logarithm of the partitioning number.

Finally, we propose the following strengthening of the Erd\H{o}s-Hajnal Conjecture that may 
potentiallly capture the real asymptotic behaviour of the EH coefficient for prime tournaments:

\begin{conjecture1}
 There exists $C>0$ such that every $n$-vertex $H$-free tournament contains a transitive subtournament 
 of order at least $n^{\frac{C}{|H|}}$.
\end{conjecture1}

As mentioned before, the Erd\H{o}s-Hajnal Conjecture is a subject of intense research however not much progress was made on it until very recently. 
In the undirected setting the Conjecture is known only for some prime graphs on at most
five vertices (\cite{safra}) and graphs obtained from them by the substitution procedure (defined in the similar way as in the
directed setting). Similarly, in the directed scenario it was known for some prime tournaments on at most five vertices
and tournaments obtained from them by the substitution procedure (see: \cite{chud2} for an excellent survey on the current state-of-the-art). 
Very recently the author of this paper together with Eli Berger and Maria Chudnovsky
proved the Conjecture for all tournaments on at most five vertices (see: \cite{bcc}). In the same paper
the Conjecture was also proven for the family of galaxies \footnote{even more recently the author of this paper proved the Conjecture for the family of so-called 
\textit{constellations} but this result has not been published yet. Furthermore, it uses similar techniques to those used in \cite{bcc} such as Szemeredi lemma; 
the methods presented in this paper can be in fact used to strengthen it}. 
The proofs used in all those previous results were of purely existential character though. The lower bounds on the EH coefficients
were very small since all those proofs relied on the Szemeredi regularity lemma. Furthermore, it was not clear at
all how to get rid of that lemma and obtain bounds that can be at least expressed by a closed-form expression.
In fact one of the big open questions was whether an expression inversely proportional to the polynomial of $|H|$
can be introduced as a universal bound for an infinite family of prime tournaments.
Some recent results (\cite{ccs}, \cite{seymour}) analyze the structure of these tournaments $H$ which exclusion implies appearance of 
the linear or almost-linear transitive subset. Those tournaments are however nonprime.
Other results focus on excluding several forbidden patterns, instead of just one. This is a much simpler scenario but even in this setting
not much is known. It is worth to mention here: \cite{seymour2}, \cite{thomasse}, \cite{zwols}.

This paper is organized as follows:
\begin{itemize}
 \item in Section 3 we define the families of stars and galaxies,
 \item in Section 4 we introduce tools used to prove all the results,
 \item in Section 5 we prove Theorem~\ref{largeenoughtheorem} and Theorem~\ref{startheorem},
 \item in Section 6 we prove Theorem~\ref{conjecturetheorem} and show how it can be used to prove
          Theorem~\ref{algorithmtheorem},
 \item in Section 7 we prove Theorem~\ref{generaltheorem},
 \item in Section 8 we prove Theorem~\ref{uppertheorem},
 \item in the Appendix we present some useful known tools for obtaining upper bounds on EH coefficients.
\end{itemize}

\section{Stars and galaxies...}

Below we define the families of stars and galaxies. All prime tournaments on at least six vertices for which the 
Conjecture is known are prime galaxies.
Stars is an important subfamily of galaxies.
The first infinite family of prime tournaments for which the Conjecture was proven were prime stars.

Fix some ordering of vertices of a tournament $H$. An edge $(v,w)$ under this ordering is called a \textit{backward edge} if $w$ precedes $v$ in this ordering.
Let T be a tournament with vertex set $V(T)$ and fix some ordering of its vertices.
The \textit{graph of backward edges}  under this ordering,
denoted by $B(T, \theta)$, has vertex set $V(T)$,  and  
$v_i v_j \in E(B(T, \theta))$ if and only if  $(v_i,v_j)$ or 
$(v_j,v_i)$ is a backward edge of $T$ under the ordering $\theta$. 
For an integer $t$, we call the graph $K_{1,t}$ a {\em star}. Let $S$ be a
star with vertex set $\{c, l_1, \ldots, l_t\}$, where $c$ is adjacent to vertices $l_1, \ldots, l_t$. We call $c$ the {\em center of the star}, and
$l_1, \ldots, l_t$ {\em the leaves of the star}. 
Note that in the case $t=1$ we may choose arbitrarily any one of the two vertices to be the center of the star, and the other vertex is then considered to be the leaf. 
Let $\theta=(v_{1},v_{2},...,v_{n})$ be an ordering of the vertex set $V(T)$ of a $n$-vertex tournament $T$. For a subset $S \subseteq V(T)$ we say that $v_{i} \in S$ is a \textit{left point of $S$ under $\theta$} if $i=\min \{j: v_{j} \in S\}$. We say that $v_{i} \in S$ is a \textit{right point of $S$ under $\theta$} if $i=\max \{j: v_{j} \in S\}$. If from the context it is clear which ordering is taken we simply say: \textit{left point of $S$} or \textit{right point of S}.  
For an ordering $\theta$ and two vertices $v_{i},v_{j}$ with $i \neq j$ we say that $v_{i}$ is \textit{before} $v_{j}$ if $i<j$ and \textit{after} $v_{j}$ otherwise.
We say that a vertex $v_{j}$ is \textit{between} two vertices $v_{i},v_{k}$ under an ordering $\theta=(v_{1},...,v_{n})$ if $i<j<k$ or $k<j<i$.

A {\em right star}
in $B(T, \theta)$ is an induced subgraph with vertex set 
$\{v_{i_0}, \ldots, v_{i_t}\}$, such that \\
$B(T,\theta)|\{v_{i_0}, \ldots, v_{i_t}\}$ is a star with center $v_{i_t}$, 
and  $i_t > i_0, \ldots, i_{t-1}$.  In this case we also  
say that $\{v_{i_0}, \ldots, v_{i_t}\}$ is  a right star in $T$.
A {\em left star}
in $B(T, \theta)$ is an induced subgraph with vertex set 
$\{v_{i_0}, \ldots, v_{i_t}\}$, such that 
$B(T,\theta)|\{v_{i_0}, \ldots, v_{i_t}\}$ is a star with center $v_{i_0}$, 
and  $i_0 < i_1, \ldots, i_t$.    In this case we also  
say that $\{v_{i_0}, \ldots, v_{i_t}\}$ is a  left star in $T$.
A {\em star} in $B(T, \theta)$ is a left star or a right star.

Let $H$ be a tournament and assume there is an ordering $\theta$ of its vertices such that every connected component of $B(H, \theta)$ is  either a star or a singleton under this ordering. 
We call this ordering a \textit{star ordering}. 

We say that a tournament is a \textit{galaxy} if there exists a star ordering of its vertices under which
no center of the star is between leaves of another star. We call such an ordering a \textit{galaxy ordering}.
If in addition  under this ordering there are no singletons then we say that a galaxy is \textit{regular}.
We say that a tournament is a \textit{star} if there exists a star ordering of its vertices under which the graph
of backward edges consists only of one connected component.

\section{Tools}

In this section we introduce several useful tools used in the consecutive parts of the paper.\\
Take a tournament $T$. 
Let $X,Y \subseteq V(T)$ be disjoint, where $|X|,|Y|>0$. Denote by $e_{X,Y}$ the number 
of directed edges $(x,y)$, where $x \in X$ and $y \in Y$.
The \textit{directed density from X to Y} is defined as 
$d(X,Y)=\frac{e_{X,Y}}{|X||Y|}.$ 
We say that a tournament $T$ is \textit{(c,$\epsilon$)-transitive} if it contains a 
transitive subtournament of order at least $c|T|^{\epsilon}$.
Define a subset $S \subseteq V(T)$ to be \textit{$c$-linear} 
if $|S| \geq c|T|$. Define a subset $S \subseteq V(T)$ to be \textit{$(c,\epsilon)$-big} 
if $|S| \geq c|T|^{\epsilon}$.  Let $\rho = (S_{1},...,S_{k})$ be a sequence of pairwise 
disjoint subsets of $V(T)$ such that for some $\lambda \geq 0$ the following holds:  
$d(S_{i},S_{j}) \geq 1 - \lambda$ for every $1 \leq i < j \leq k$. We say that $\rho$ 
is an \textit{$(c,\lambda)$-$l$-sequence} if every $S_{i}$ is $c$-linear. We say $\rho$ 
is an \textit{$(c,\lambda,\epsilon)$-$t$-sequence} if every $S_{i}$ is transitive and 
$(c,\epsilon)$-big. We say that $\rho$ is an \textit{$(c_{1},c_{2},\lambda,\epsilon)$-$m$-sequence} 
if: $S_{i}$ is $c_{2}$-linear for odd $i$ and $S_{i}$ is transitive and $(c_{1},\epsilon)$-big for 
even $i$.
We say that $\rho$ is smooth if for every $1 \leq i < j \leq k$ we have: 
$d(\{v\},S_{j}) \geq 1 - \lambda$ for all $v \in S_{i}$ and $d(S_{i},\{v\}) \geq 1 - \lambda$ 
for all $v \in S_{j}$.



Whenever we do not care about parameters of the $(c,\lambda,\epsilon)$-$t$-sequences, 
$(c,\lambda)$-$l$-sequences or $(c_{1},c_{2},\lambda,\epsilon)$-$m$-sequences under consideration, 
we refer to them simply as: $t$-sequences, $l$-sequences and $m$-sequences respectively.

\section{Proof of Theorem~\ref{largeenoughtheorem} and Theorem~\ref{startheorem}}

We begin by proving Theorem~\ref{largeenoughtheorem}.\\

\Proof
Assume that the lemma is false. Then there exists an $H$-free tournament $T_{H}$ such that its largest transitive
subtournament is of order smaller than $|T_{H}|^{\epsilon(H)}$.
But then, following the proof of Theorem \ref{strongerhomotheorem} from Appendix A, we conclude that
there exists an infite family of $H$-free tournaments: $\{F_{0}, F_{1}, ...\}$ such that each $F_{i}$
does not contain transitive subtournaments of order larger than 
$|F_{i}|^{\frac{\log(tr(T_{H}))}{\log(T_{H})}}$. 
Since $tr(T_{H}) < |T_{H}|^{\epsilon(H)}$ we conclude that there exists $\epsilon > 0$ such that
$tr(|F_{i}|) \leq |F_{i}|^{\epsilon(H) - \epsilon}$ for all $F_{i}$'s. That clearly contradicts the fact that
every $n$-vertex $H$-free tournament (in particular, every $F_{i}$) contains a transitive subtournament of order at least $c(H)n^{\epsilon(H)}$.
\bbox

Now we prove  Theorem~\ref{startheorem}.

\Proof
Let $H$ be a prime star. Note that we know that $\epsilon(H) = O(\frac{\log(|H|)}{|H|})$ (see: Appendix A).
We will prove that every $H$-free tournament contains a transitive subtournament
of order at least $c(H) n^{\epsilon(H)}$ for $\epsilon(H)=\frac{1}{3|H|\log(2|H|)}$ and $c(H) = \min(4^{-|H|}, 2^{-\frac{1}{1-\epsilon(H)}})$.
That, according to Theorem~\ref{largeenoughtheorem}, completes the proof.
We proceed by induction on the order of the $H$-free tournament $T$. 
The statement is trivial for tournaments $T$ with no more than $\max(4^{|H|}, 2^{\frac{1}{1-\epsilon(H)}})$ vertices.
Now let $T$ be a $n$-vertex $H$-free tournament for $n>\max(4^{|H|},2^{\frac{1}{1-\epsilon(H)}}))$. 
Note that by Ramsey thery, $T$ contains
a transitive subtournament of order at least $2|H|$.
Let $h_{c}$ be a center of the star $H$ and let $(h_{1},...,h_{l})$ be a transitive ordering of 
its remaining vertices. Denote by $L$ the largest 
transitive subtournament of $T$
and let $W = V(T) \backslash V(L)$. We have: $|L| \geq 2l$.
If $|W| = 0$ then we are done. Thus assume that $|W|>0$.
Partition $V(L)$ into $l+1$ transitive subsets: $L_{1},...,L_{l}, R$, 
each $L_{i}$ of size $\lfloor \frac{|L|}{l} \rfloor$,
such that for every $1 \leq i < j \leq l$ $L_{i}$ is outadjacent to  $L_{j}$
and each $L_{i}$ is adjacent to $R$.
Consider the following mapping: $\zeta: W \times \{1,...,l\} \rightarrow \{0,1\}$: 

$$\zeta(w,i) = \begin{cases}
  1 & \text{if $(h_{c},h_{i})$ and $\exists_{r \in L_{i}}$ such that $(w,r)$,} \\
  1 & \text{if $(h_{i},h_{c})$ and $\exists_{r \in L_{i}}$ such that $(r,w)$,} \\  
  0 & \text{otherwise.} 
\end{cases}
$$
Note that if there exists $w \in W$ such that $\zeta(w,i) = 1$ for $i=1,...,l$ then $T$ contains the copy of $H$.
Therefore we can assume that for every $w \in W$ there exists $i$ such that $\zeta(w,i) = 0$. But then, by Pigeonhole
Principle, there exists $i_{0} \in \{1,...l\}$ and a subset $S \subseteq W$ of size at least $\frac{|W|}{l}$ such that 
for every $w \in S$ we have: $\zeta(w,i_{0})=0$. Notice that this implies that either every vertex of $S$ is adjacent
to every vertex of $L_{i_{0}}$ or every vertex of $L_{i_{0}}$ is adjacent to every vertex of $S$.
Let $F$ be the largest transitive subtournament of a tournament induced by $S$. By induction we have:
$|F| \geq c(H)|S|^{\epsilon(H)}$. Now note that the tournament obtained by merging $F$ with $L_{i_{0}}$
is transitive. Its size is: $|F| + \lfloor \frac{|L|}{l} \rfloor$. Thus, from the definition of
$L$ we get: $|F| + \lfloor \frac{|L|}{l} \rfloor \leq |L|$. So:

\begin{equation}
\label{eq0}
c(H)(\frac{n-|L|}{l})^{\epsilon(H)} + \lfloor \frac{|L|}{l} \rfloor \leq |L|.
\end{equation}

We can assume that $|L| \leq \frac{n}{2}$ since otherwise, the condition $n > 2^{\frac{1}{1-\epsilon(H)}}$
gives us: $|L| \geq n^{\epsilon(H)}$ and we are done.
We also have: $\lfloor \frac{|L|}{l} \rfloor \geq \frac{|L|}{l} - 1 \geq  \frac{|L|}{2l}$, where  
the last inequality comes from the fact that $|L| \geq 2l$.
Thus from \ref{eq0} we get: $|L| \geq c(H)\frac{1}{1-\frac{1}{2l}}\frac{n^{\epsilon(H)}}{(2l)^{\epsilon(H)}}$.
To complete the proof it suffices to notice that under our choice of $\epsilon(H)$ we have:
$(1-\frac{1}{2l}) (2l)^{\epsilon(H)} \leq 1$.

\bbox

\section{Proof of Theorem ~\ref{conjecturetheorem} and Theorem~\ref{algorithmtheorem}}

In this section we prove the main theorem of the paper, Theorem~\ref{conjecturetheorem}
and show how Theorem~\ref{algorithmtheorem} can be derived from it. We start with 
Theorem~\ref{conjecturetheorem}.\\

\Proof 
Let $H$ is a prime galaxy. It is easy to see that without loss of generality 
we can assume that $H$ is a regular galaxy.
Denote $V(H)=\{h_{1},...,h_{|H|}\}$. 
Take a galaxy ordering of the vertices of $H$. Denote this ordering by $(h_{1},...,h_{|H|})$.
Assume that first $k_{1}$ vertices of this ordering are centers of stars, next
$w_{1}$ are leaves, next $k_{2}$ are centers of stars, next $w_{2}$ are leaves, ...
and finally last $k_{t+1}$ are centers of stars, where: $t$ is some nonnegative integer.
Denote by $g$ the number of stars of $H$ (note that we have: $k_{1} + ... + k_{t+1}=g$)
and let $r = h - g$.
Let $T$ be an $n$-vertex $H$-free tournament. 

We will proceed by induction on $n$.
Note first that for every given constant $C>0$ ($C$ does not depend on $H$ and $n$) we can assume that
the theorem holds for $n \leq C \cdot 2^{|H|^{4}}$ by taking the constant hidden in  
the $O(\frac{1}{|H|^{5}})$ expression to be small enough.
Then from now one we will assume that $n > C \cdot 2^{|H|^{4}}$, where $C$ is taken to be large enough
(but does not depend on $n$ and $H$).
We use standard notation: $h=|H|$.
We need to prove first the following lemma:

\begin{lemma}
 \label{msequencelemma}
 Assume that an $H$-free tournament $T$ contains an
 $m$-sequence $\rho=(S_{1},T_{1},...,T_{t},S_{t+1})$ which is smooth. Assume that $\rho$ is an 
 $(c_{1},c_{2},\lambda,\epsilon)$-$m$-sequence, for $\lambda \leq \frac{1}{24h^{3}2^{(t+1)h^{2}+h}}$, 
 some parameters $c_{1}, c_{2}$,  
 and $\epsilon= \frac{\log(1+\frac{1}{M_{2}-1})}{\log(M_{1})}$, where:
 $M_{1}=\frac{24 \cdot 2^{k_{1}^{2}+...+k_{t+1}^{2}+g}}{c_{2}}$ and
 $M_{2}=\frac{4 \max_{i} w_{i}}{c_{1}}$.
 Assume furthermore that: $n \geq \frac{2h}{c_{1}}$ and $(1-\frac{eh}{c_{2}n})^{h} \geq \frac{1}{2}$.
 Then the following holds: if every proper subtournament of $T$ is 
 $(1,\epsilon)$-transitive, $T$ is also $(1,\epsilon)$-transitive.
\end{lemma}

\Proof
Note first that every $n$-vertex tournament has at least $\delta(k)=\eta\frac{1}{2^{k^{2}}}n^{k}$ transitive 
subtournaments of order $k$, where: $\eta = (1-\frac{k}{n})^{k}$.
Indeed, the number of transitive subtournaments is at least $\frac{{n \choose 2^{k}}}{{n-k \choose 2^{k}-k}}$.
By evaluating this expression we get the formula for $\delta(k)$.
Now we need some generalization of this result.
For given $l_{1},...,l_{j}>0$, $t_{1},...,t_{j}>0$, $\lambda>0$ denote by 
$\theta(t_{1},...,t_{j},l_{1},...,l_{j},\lambda)$
the number such that every smooth $(c,\lambda)$-$l$-sequence $(L_{1},...,L_{j})$ with
$|L_{1}| = l_{1}$,...,$|L_{j}|=l_{j}$ contains at least 
$\theta(t_{1},...,t_{j},l_{1},...,l_{j},\lambda)$ transitive subtournaments 
of order $t_{1} + ... + t_{j}$ such that
their first $t_{1}$ vertices under transitive ordering are in $L_{1}$, next $t_{2}$ are in 
$L_{2}$, etc.
Let $Tr_{1}$ be some transitive tournament found in the tournament induced by $L_{1}$.
Denote by $W_{i}$ for $i=2,...,j$ the subset of $L_{i}$ that consists of vertices adjacent 
from all the vertices of $Tr_{1}$. Note that $|W_{i}| \geq l_{i}(1-t_{1}\lambda)$.
Note also that if we find in the $l$-sequence $(W_{2},...,W_{j})$ a transitive 
tournament $Tr_{2}$
of order $t_{2}+...+t_{j}$ with first $t_{2}$ vertices under transitive ordering in $W_{2}$,
next $t_{3}$ in $W_{3}$, etc. Then by merging $T_{1}$ with $T_{2}$ we obtain a tournament
that is counted by $\theta(t_{1},...,t_{j},l_{1},...,l_{j},\lambda)$.
From our previous remark we know that the number of tournaments $Tr_{1}$ is at least
$\eta_{1}\frac{1}{2^{t_{1}^{2}}}l_{1}^{t_{1}}$, 
where $\eta_{1} = (1-\frac{t_{1}}{l_{1}})^{t_{1}}$.
Thus we get the following simple recurence formula:
$$\theta(t_{1},...,t_{j},l_{1},...,l_{j},\lambda) \geq 
\eta_{1}\frac{1}{2^{t_{1}^{2}}}l_{1}^{t_{1}} \theta(t_{2},...,t_{j},l_{1}(1-t_{1}\lambda),...,l_{j}(1-t_{1}\lambda),\frac{\lambda}{1-t_{1}\lambda}).$$
For fixed $j$ and $t_{1},...,t_{j}$, if we take 
$l_{1},...,l_{j}$ large enough and solve the recurence above, we get:
$$\theta(t_{1},...,t_{j},l_{1},...,l_{j},\lambda) \geq 
\frac{l_{1}^{t_{1}} \cdot...\cdot l_{j}^{t_{j}}}{2 \cdot 2^{t_{1}^{2}+...+t_{j}^{2}}}
((1-t_{1}\lambda_{max}) \cdot ... \cdot (1-t_{j}\lambda_{max}))^{t_{1}+...+t_{j}},$$ where: 
$\lambda_{max}$ satisfies: $\lambda_{max \geq }\frac{\lambda}{(1-t_{1}\lambda_{max})...(1-t_{j}\lambda_{max})}$.
All $l_{i}$'s should be large enough to satisfy: $(1-\frac{h}{p_{i}})^{h} \geq \frac{1}{2}$, 
where: $p_{i}=l_{i}((1-t_{1}\lambda_{max})...(1-t_{j}\lambda_{max}))^{k_{1}+...+k_{j}}$.

If we assume besides that $\lambda_{max} \leq \frac{1}{2\max_{i} t_{i}}$, then using the inequality
$1 - x \geq e^{-2x}$ (for $x \leq \frac{1}{2}$) we get:

$$\theta(t_{1},...,t_{j},l_{1},...,l_{j},\lambda) \geq \frac{l_{1}^{t_{1}}...l_{j}^{t_{j}}}{2 \cdot 
2^{t_{1}^{2}+...+t_{j}^{2}}} (e^{-2t_{1}\lambda_{max}-...-2t_{j}\lambda_{max}})^{g}.$$

Thus we get:

$$\theta(t_{1},...,t_{j},l_{1},...,l_{j},\lambda) \geq \frac{l_{1}^{t_{1}}...l_{j}^{t_{j}}}{2 \cdot 
2^{t_{1}^{2}+...+t_{j}^{2}}} e^{-2g^{2}\lambda_{max}}.$$

Now, if we take sets: $S_{1},...,S_{t+1}$ as $L_{1},...,L_{j}$, the sequence: $k_{1},...,k_{t+1}$ as $t_{1},...,t_{j}$
and denote $\theta = \theta(|S_{1}|,...,|S_{t+1}|,k_{1},...,k_{t+1},\lambda)$
we obtain (under previous assumptions):

$$\theta \geq \alpha |S_{1}|^{k_{1}}...|S_{t+1}|^{k_{t+1}},$$ where:
$\alpha = \frac{e^{-2g^{2}\lambda_{max}}}{2^{k_{1}^{2}+...+k_{t+1}^{2}+1}}$.

Note that we have to assume that $|S_{1}|,...,|S_{k}|$ are large enough.
Lets see how large. According to previous remarks we need:
$(1-\frac{h}{P_{i}})^{h} \geq \frac{1}{2}$, where:
$P_{i} = |S_{i}|((1-k_{1}\lambda_{max})...(1-k_{t+1}\lambda_{max}))^{k_{1}+...+k_{t+1}}$.
We also need:
\begin{itemize}
 \item $\lambda_{max} \leq \frac{1}{2\max_{i \in \{1,...,t+1\}}k_{i}}$,
 \item $\lambda \leq \lambda_{max} (1-k_{1}\lambda_{max})...(1-k_{t+1}\lambda_{max})$.
\end{itemize}

Under given assumption on $\lambda_{max}$ the first assumption can be replaced by:
$(1-\frac{h}{P})^{h} \geq \frac{1}{2}$, where: $P=c_{2}n e^{-2g^{2}\lambda_{max}}$.
We call this condition the \textit{strong linearity condition} since it says that the
linear sets we start with in the $m$-sequence must be large enough. Assume that this condition holds
(we will see later why it is true under lemma assumptions).

Now we will try to construct $H$ in our $(c_{1},c_{2},\lambda,\epsilon)$-$m$-sequence star by star and show
that if we cannot succed then we get big enough transitive subtournament. We can divide each transitive
chunk $T_{i}$ of our $m$-sequence into $w_{i}$ subchunks of the same size 
$\lfloor \frac{|T_{i}|}{w_{i}} \rfloor$ (and get rid of its last 
$|T_{i}|-w_{i}\lfloor \frac{|T_{i}|}{w_{i}} \rfloor$ vertices under transitive ordering) in such a way that the first one consists
of first $\lfloor \frac{|T_{i}|}{w_{i}} \rfloor$ vertices of $T_{i}$ under its transitive ordering, next one consists
of next $\lfloor \frac{|T_{i}|}{w_{i}} \rfloor$ of its vertices under its transitive ordering and so on.
Note that $q=\lfloor \frac{|T_{i}|}{w_{i}} \rfloor \geq \frac{|T_{i}|}{w_{i}}-1 \geq 
\frac{|T_{i}|}{2w_{i}}$, where the last inequality is true if $|T_{i}| \geq 2h$.
Thus we need: $c_{1}n \geq 2h$, i.e. $n \geq \frac{2h}{c_{1}}$, but this one of the assumptions
of the lemma. 
Let us order the stars of $H$ as follows: $\Sigma_{1},...,\Sigma_{g}$ (we can assume without loss of generality
that $H$ has no singletons). 
Notice that the $i^{th}$ subchunk corresponds to the $i^{th}$ leaf under given galaxy ordering.
Notice also that the set of all the centers of stars of $H$ is a transitive set.
Denote by $\mathcal{T}$ the set of all the transitive tournaments of $k_{1}+...+k_{t+1}$ vertices each
and such that the first $k_{1}$ vertices under transitive ordering of each of them are in $S_{1}$, next
$k_{2}$ are in $S_{2}$, etc. We have already proved that the number of all of them is at least 
$\alpha|S_{1}|^{k_{1}}...|S_{t+1}|^{k_{t+1}}$.
We associate with the $i^{th}$ center for $i=1,...,g$ the set of $i^{th}$ vertices of tournaments from 
$\mathcal{T}$ under their transitive ordering. 
Let us take the first star $\Sigma_{1}$. Without loss of generality assume that its center is  the first center 
of the given galaxy ordering and the star is a left star. Let $\mathcal{T}^{1}_{\frac{\alpha}{2}}$ be the set of those vertices $v$ of tournaments
from $\mathcal{T}$ that are associated with the center of $\Sigma_{1}$ and such that for each of them there are at
least $\frac{\alpha}{2}|S_{1}|^{k_{1}-1}|S_{2}|^{k_{2}}...|S_{t+1}|^{k_{t+1}}$ tournaments of $\mathcal{T}$
with $v$ being their first vertex under transitive ordering. Simple counting argument give us:
$|\mathcal{T}^{1}_{\frac{\alpha}{2}}| \geq \frac{\alpha}{2-\alpha}|S_{1}|$.
For any vertex $v \in \mathcal{T}^{1}_{\frac{\alpha}{2}}$ consider subchunks of all the sets $T_{i}$ ($i=1,...,t$) 
that are associated with leaves of $\Sigma_{1}$. Denote those subchunks as: $ST^{v}_{1},...,ST^{v}_{w_{1}}$.
If for some $v \in \mathcal{T}^{1}_{\frac{\alpha}{2}}$ in every $ST^{v}_{i}$ there exists a vertex $y_{i}$ adjacent to $v$ then the set $\{v,y_{1},...,y_{w_{1}}\}$
induces $\Sigma_{1}$. We call this case: "\textit{the star setting}". We then delete all the transitive subchunks related to the leaves of $\Sigma_{1}$.
We also modify other transitive subchunks in the following way.
Let $ST_{*}$ be one of the other chunks. For $x \in \{v, y_{1},..., y_{w_{1}}\}$ we denote by $N^{x}_{ST_{*}}$:
\begin{itemize}
\item the set of vertices of $ST_{*}$ adjacent from $x$ if the vertex of $\Sigma_{1}$ that $x$ corresponds to is before 
the leaf of $\Sigma_{1}$ that $ST_{*}$ corresponds to 
\item the set of vertices of $ST_{*}$ adjacent to $x$ otherwise.
\end{itemize}

Let us denote $u(ST_{*}) = \bigcap_{x \in \{v,y_{1},...,y_{w_{1}}\}} N^{x}_{ST_{*}}$.
We replace each $ST_{*}$ by $u(ST_{*})$. From the definition of 
($c_{1}$,$c_{2}$,$\lambda$, $\epsilon$)-$m$-sequencewe have:
$|u(ST_{*})| \geq |ST_{*}| - h_{1}\lambda T_{i^{*}}$, where $h_{1}$ is the number of vertices of $\Sigma_{1}$
and $T_{i^{*}}$ is the transitive element of the given $m$-sequence that $ST_{*}$ belongs to. 
Denote by $\mathcal{T}^{'}$ the set of tournaments from $\mathcal{T}$
such that their first vertex under transitive ordering is $v$. We have already showed that:
$|\mathcal{T}^{'}| \geq \frac{\alpha}{2}|S_{1}|^{k_{1}-1}|S_{2}|^{k_{2}}...|S_{t+1}|^{k_{t+1}}$.
Denote by $\mathcal{T}^{''}$ the subset of $\mathcal{T}^{'}$ consisting of those transitive tournaments 
of $\mathcal{T}^{'}$ such that for each of them there exists a vertex $x \in \{y_{1},...,y_{w_{1}}\}$ with 
the following property:
\begin{itemize}
 \item $x$ is adjacent to some vertex $w \in V(\mathcal{T}^{'})$ and belongs to the subchunk
       appearing later in the $m$-sequence than a linear set from which $w$ was taken, or
 \item $x$ is adjacent from some vertex $w \in V(\mathcal{T}^{'})$ and belongs to the subchunk
       appearing earlier in the $m$-sequence than a linear set from which $w$ was taken. 
\end{itemize}

From the definition of the $(c_{1},c_{2},\lambda,\epsilon)$ $m$-sequence 
we get: $$|\mathcal{T}^{''}| \leq \sum_{i=1}^{t+1} |S_{1}|^{k_{1}-1}|S_{2}|^{k_{2}}...|S_{t+1}|^{k_{t+1}} 
(\frac{h_{1}h\lambda|S_{1}|}{|S_{1}|}+...+\frac{h_{1}h\lambda|S_{t+1}|}{|S_{t+1}|}).$$
Thus we get $|\mathcal{T} \backslash \mathcal{T}^{''}| \geq 
(\frac{\alpha}{2}-(t+1)h_{1}h\lambda)|S_{1}|^{k_{1}-1}|S_{2}|^{k_{2}}...|S_{t+1}|^{k_{t+1}}$.
We replace $\mathcal{T}$ by $\mathcal{T} \backslash \mathcal{T^{''}}$ and replace all the chunks $ST_{*}$ 
that were not already  removed with the leaves of the star $\Sigma_{1}$ by $u(ST_{*})$.
We then proceed in the analogous way for the star $\Sigma_{2}$. 
On the other hand, if we do not have a "star setting" then, by Pigeonhole principle, we know that at least 
$\frac{1}{h_{1}}|\mathcal{T}^{1}_{\frac{\alpha}{2}}|$ of the vertices of $\mathcal{T}^{1}_{\frac{\alpha}{2}}$
are complete to/from some subchunk of the transitive set in the given $m$-sequence. We call this setting
a "\textit{non-star setting}".
If we encounter a "star setting" every time we are looking for the star then we can merge all the stars that were found by us so far. The way we update the
entire $m$-sequence enables us to conclude that by merging all those stars we get a copy of $H$, a contradiction.
Thus at some point we get a "non-star setting". From our earlier analysis it is clear that if this is the case
then we get a set of size at least 
$\frac{\alpha_{f}}{2-\alpha_{f}}\frac{1}{h}c_{2}n$ complete to/from a transitive set of size at least
$(\frac{1}{2\max_{i}w_{i}}-h_{1}\lambda - h_{2} \lambda -... - h_{g}\lambda)c_{1}n^{\epsilon}$, where: $\alpha_{f} = (\frac{\alpha}{2^{g}}-(t+1)h^{2}\lambda)$ and $h_{i}$'s are sizes of stars.
Denote $w = \max_{i} w_{i}$, $A = \frac{1}{w} - h\lambda$, 
$B=\frac{\alpha_{f}}{(2-\alpha_{f})h}c_{2}$. 
We conclude that we got a transitive tournament of order
at least $Ac_{1}n^{\epsilon}$ complete from/to the linear set of size at least $Bn$. By the assumptions of the lemma,
we know that this linear set contains a transitive tournament $R$ of order at least $(Bn)^{\epsilon}$.

If we merge it with a tournament of order at least $Ac_{1}n^{\epsilon}$, then we get a transitive tournament of order at 
least $Ac_{1}n^{\epsilon}+(Bn)^{\epsilon}$. Thus we get: 
$tr(T) \geq Ac_{1}n^{\epsilon}+(Bn)^{\epsilon}$, i.e.: $tr(T) \geq (Ac_{1} + B^{\epsilon})n^{\epsilon}$.
Take $K_{1} = B^{-1}$, $K_{2}=A^{-1}$. It is easy to see that for 
$\epsilon \leq \frac{\log(1 + \frac{1}{K_{2}-1})}{\log(K_{1})}$ we have: 
$Ac_{1}+B^{\epsilon} \geq 1$. Thus for such a choice of $\epsilon$ 
we get: $tr(T) \geq n^{\epsilon}$.
Let us summarize our assumptions. We have: $\alpha = \frac{e^{-2g^{2}\lambda_{max}}}
{2^{k_{1}^{2}+...+k_{t+1}^{2}+1}}$, $K_{1}=\frac{2-\frac{\alpha}{2^{g}}+
\lambda(t+1)h^{2}}{c_{2}(\frac{\alpha}{2^{g}}-\lambda(t+1)h^{2})}$,
$K_{2}=\frac{1}{(\frac{1}{2w}-h\lambda)c_{1}}$. We also need to assume that:
$\lambda \leq \frac{1}{2\max k_{i}}$ and $\lambda \leq 
\lambda_{max}(1-k_{1}\lambda_{max})...(1-k_{t+1}\lambda_{max})$. Note that $1 - k_{i} \lambda_{max}
 \geq e^{-2k_{i}\lambda_{max}}$. Thus it suffices to have:
 $\lambda \leq \lambda_{max} e^{-2k_{1} \lambda_{max}-...-2k_{t+1}\lambda_{max}}$ and
 $\lambda_{max} \leq \frac{1}{2\max k_{i}}$. Therefore it is enough to have: 
 $\lambda \leq \lambda_{max} e^{-2g\lambda_{max}}$ and $\lambda_{max} \leq \frac{1}{2 \max k_{i}}$.
 We also want the following inequality: $\lambda h^{2}(t+1) \leq \frac{\alpha}{2^{g+1}}$, i.e. 
 $\lambda \leq \frac{\alpha}{2^{g+1}h^{2}(t+1)}$. For this choice of $\lambda$ and $\lambda_{max}$ we get:
 $K_{1} \leq \frac{4}{c_{2}\frac{\alpha}{2^{g}}}$. If we furthermore have: 
 $\lambda_{max} \leq \frac{1}{2g^{2}}$, then we obtain: $\alpha \geq 
 \frac{e^{-1}}{2^{k_{1}^{2}+...+k_{t+1}^{2}+1}}$. Then we know that: 
 $K_{1} \leq \frac{8e \cdot 2^{k_{1}^{2}+...+k_{t+1}^{2}+g}}{c_{2}}$.
 Thus: $K_{1} \leq \frac{24 \cdot 2^{k_{1}^{2}+...+k_{t+1}^{2}+g}}{c_{2}}$.
 For $h\lambda \leq \frac{1}{4w}$ we also obtain: $K_{2} \leq \frac{4w}{c_{1}}$.
 Taking into account all the inequalities on $\lambda$ and $\lambda_{max}$ we derived so far, it is easy to
 see that those upper bounds on $K_{1}$ and $K_{2}$ are valid for 
 $\lambda \leq \frac{1}{24h^{3}2^{(t+1)h^{2}+h}}$ (we leave this simple check to the reader). 
 It is also easy to check that for $\lambda_{max} \leq \frac{1}{2g^{2}}$ and under lemma assumptions
 the strong linearity condition holds.
 That completes the proof of Lemma~\ref{msequencelemma}.
 \bbox
Now we state and prove another useful lemma:

\begin{lemma}
\label{lintotrans}
 Assume that $T$ contains a smooth $(c,\lambda_{0})$-$l$-sequence $\chi=(L_{1},...,L_{2t+1})$. 
 Assume that every proper subtournament of $T$ is $(1,\epsilon)$-transitive, where:
 $\epsilon \leq \frac{\log(2)}{\log(\frac{2}{c})}$.
 Then $T$ contains a smooth $(\frac{1}{4},\frac{c}{2},\lambda,\epsilon)$-$m$-sequence of 
 length $2t+1$, where: $\lambda=64t^{2}(t+1)\lambda_{0}$.
\end{lemma}

\Proof

In this proof we will very often use terms: "transitive tournament" and "transitive set" 
interchangeably since the context will be always clear.
Take some $L_{2i}$ for $i=1,2,...,t$. Since a tournament induced by $L_{2i}$ is 
$(1,\epsilon)$-transitive,
it contains a transitive subtournament $T^{i}_{1}$ or 
order $\lceil (\frac{cn}{2})^{\epsilon} \rceil$. We delete $T^{i}_{1}$ and repeat the 
procedure to get $T^{i}_{2}$. We continue as long as the size of the set of vertices
remaining in $L_{2i}$ is at least $\frac{|L_{2i}|}{2}$. 
Denote by $L_{2i}^{'}$ the set of vertices obtained from merging deleted transitive tournaments
$T^{i}_{j}$ for $j=1,...$. Note that under our assumption on $\epsilon$ we have:
$|T^{i}_{j}| \geq \frac{1}{2}n^{\epsilon}$.
Now notice that  a sequence obtained from $\chi$ by replacing every $L_{2i}$ by $L_{2i}^{'}$
for $i=1,2,...,t$
is a smooth $(\frac{c}{2},2\lambda_{0})$-$m$-sequence (simple density argument).
Denote this sequence as $\chi^{'}$.
Denote by $t_{i}$ for $i=1,2,...,t$ the number of transitive tournaments $T^{i}_{j}$ creating
$L_{2i}^{'}$. Denote by $n_{i}$ those tournaments $T^{i}_{j}$ that satisfy the following:
there exists a set $L_{j}$ such that
\begin{itemize}
 \item $j < i$ and $d(L_{j},L_{i}) < 1 - 2 \lambda_{0} W$ or
 \item $j > i$ and $d(L_{i},L_{j}) < 1 - 2 \lambda_{0} W$,
\end{itemize}

where $W=4t$. 
We call tournaments $T^{i}_{j}$ with this property $W$-\textit{bad}.
Tournaments $T^{i}_{j}$ that do not have this property will be denoted as $W$-\textit{good}.
Again, a simple density argument gives us: $n_{i} \leq 2t \frac{t_{i}}{W}$.
From our choice of $W$ we get: $n_{i} \leq \frac{t_{i}}{2}$. Replace in the $m$-sequence 
$\chi^{'}$ every $L_{2i}^{'}$ by its subset obtained by taking all related
$W$-good transitive tournaments (i.e. by getting rid of $W$-bad transitive tournaments).
Denote the new $m$-sequence constructed in such a way as $\chi^{''}$. Note that $\chi^{''}$
is a smooth $(\frac{c}{2},8t\lambda_{0})$-$l$-sequence.
If we can now find $W$-good transitive tournaments: $T^{1}_{j_{1}},...,T^{t}_{j_{t}}$ satisfying
for $i_{1} < i_{2}$:
$d(T^{i_{1}}_{j_{i_{1}}},T^{i_{2}}_{j_{i_{2}}}) \geq 1 - 8Mt\lambda_{0}$ (for some constant $M>0$)
then $(L_{1},T^{1}_{j_{1}},L_{2},...,T^{t}_{j_{t}},S_{t+1})$ is a 
$(\frac{1}{2},c,\lambda_{1},\epsilon)$-$m$-sequence for $\lambda_{1}=8Mt\lambda_{0}$.
Let us construct the $k$-partite graph $G$ with color classes: $A_{1},...,A_{t}$ such that the
vertices of $A_{i}$ are $W$-good transitive tournaments $T^{i}_{j}$ and there exists an edge
between $x \in A_{i}$ and $y \in A_{j}$ for $i<j$ iff $d(V(x),V(y)) \geq 1 - 8Mt\lambda_{0}$.
From the simple density argument we know that in this graph there are always at least
$(1-\frac{1}{M})|A_{i}||A_{j}|$ edges between any vertex of $A_{i}$ and a set $A_{j}$. Note that a clique
of order $t$ in this graph corresponds to the sequence $T^{1}_{j_{1}},...,T^{t}_{j_{t}}$.
To construct a clique in $G$ of size $t$ we choose an arbitrary vertex $v_{1}$ in $A_{1}$ and replace
$A_{2},...,A_{t}$ by the sets of its neighbors. We then choose an arbitrary vertex in the
set of neighbors of $v_{1}$ in $A_{2}$ and repeat the entire procedure. It is easy to see that
we will succeed if $M>t$. Thus we can conclude, using our previous remarks, that for 
$\lambda_{1} =  8t(t+1)\lambda_{0}$ we obtain an $(\frac{1}{2},c,\lambda_{1},\epsilon)$-$m$-sequence.
Denote it as: $(F_{1},...,F_{2t+1})$. For every $F_{i}$ denote by $F^{b}_{i}$ the subset
of vertices of $v$ of $F_{i}$ that satisfy the following: there exists $j$ such that
\begin{itemize}
 \item $j<i$ and $d(F_{j},\{v\}) < 1 - W\lambda_{1}$ or
 \item $j>i$ and $d(\{v\},F_{j}) < 1 - W\lambda_{1}$.
\end{itemize}

As previously, we conclude (using simple density analysis) that:
$|F^{b}_{i}| \leq 2t \frac{|F_{i}|}{W}$. Denote $Q_{i} = F_{i} \backslash F^{b}_{i}$.
We get $|Q_{i}| \geq \frac{1}{2}|F_{i}|$. Not it is easy to see that 
$(Q_{1},...,Q_{2t+1})$ is a smooth $(\frac{1}{4},\frac{c}{2},\lambda,\epsilon)$-$m$-sequence,
where: $\lambda=2W\lambda_{1}=64t^{2}(t+1)\lambda_{0}$. That completes the proof.
\bbox

We will now introduce an important parameter $C^{h}(c,\lambda)$ having the following property:

for every set of subsets $\{S_{1},...,S_{h}\}$ of the set $V(Z)$, where $Z$ is some tournament,
if $|S_{i}| \geq c |Z|$ for $i=1,...,h$ the the following holds:
\begin{itemize}
 \item there exists a sequence $s_{1},...,s_{h}$ s.t. $s_{i} \in S_{i}$ and there exists a mapping:
       $\phi: s_{i} \rightarrow h_{i}$ which is an isomorphism between a tournament induced by 
       $\{s_{1},...,s_{h}\}$ and $H$ or
 \item there exist in $V(Z)$ two disjoint subsets: $X_{1}$ and $X_{2}$ such that 
       $|X_{1}|,|X_{2}| \geq C^{h}(c,\lambda)|Z|$ and $d(X_{1},X_{2}) \geq 1 - \lambda$.
\end{itemize}

\begin{lemma}
 \label{shortlemma}
 We can take $C^{h}(c,\lambda) = \frac{\lambda^{h}c}{h}$.
\end{lemma}

\Proof
Choose an arbitrary vertex $s_{1} \in S_{1}$. Denote by $N^{+}_{s_{1},S_{j}}$ for $j \neq 1$ the
number of outneighbors of a vertex $s_{1}$ in $S_{j}$ and by $N^{-}_{s_{1},S_{j}}$ the number of
inneighbors of a vertex $s_{1}$ in $S_{j}$.
Assume first that for every $j \neq 1$ we have: $|N^{+}_{s_{1},S_{j}}| \geq \lambda |S_{j}|$
and $|N^{-}_{s_{1},S_{j}}| \geq \lambda |S_{j}|$. We call this setting \textit{a regular setting}.
If the regular setting holds then for $j \neq 1$ we define $B_{s_{1},S_{j}}$ to be 
$N^{+}_{s_{1},S_{j}}$ if $h_{1}$ is adjacent to $h_{j}$ and to be $N^{-}_{s_{1},S_{j}}$ otherwise.
Now note that if one can find vertices $s_{2},...,s_{h}$ such that $s_{j} \in B_{s_{1},S_{j}}$
for $j=2,...,h$ with the property that there exists an isomorphism $\psi: s_{j} \rightarrow h_{j}$ for
$j=2,...,h$ then we notice that a set $\{s_{1},...,s_{h}\}$ induces a copy of $H$. Furthermore, the isomorphism
$\phi$ between a tournament induced by $\{s_{1},...,s_{h}\}$ and $H$ is defined by the mapping $s_{j} \rightarrow h_{j}$ for $j=1,...,h$.
If no vertices $s_{2},...,s_{h}$ with this property can be found then by the definition of $C$ we get
two disjoint sets of vertices $D_{1}, D_{2}$ such that 
$|D_{1}|,|D_{2}| \geq C^{h-1}(\lambda c, \lambda)|Z|$ and $d(D_{1},D_{2}) \geq 1 - \lambda$.
Now assume that we do not have a regular setting. Then, by the Pigeonhole principle, there exists
$j_{*} \neq 1$ such that at least $\frac{c|Z|}{h}$ vertices $v$ of $S_{1}$ satisfy:
$d(S_{j_{*}},\{v\}) \geq 1 - \lambda$ if $h_{1}$ is adjacent to $h_{j_{*}}$ and
$d(\{v\},S_{j_{*}}) \geq 1 - \lambda$ if $h_{j_{*}}$ is adjacent to $h_{1}$.
Thus we get the following recurrence: $C^{h}(c,\lambda) = 
\min(\frac{c}{h},C^{h-1}(\lambda c, \lambda))$ for $h \geq 1$. We can also obviously assume that:
$C^{0}(c,\lambda)=1$. solving this recurrence gives us: $C^{h}(c,\lambda) = \frac{\lambda^{h}c}{h}$.

\bbox

To finalize the proof of Theorem~\ref{conjecturetheorem} we need one more technical lemma.

\begin{lemma}
\label{finallemma}
For any $u>0$ and $0 < \lambda <1$ the following holds:
if a tournament $T$ is $H$-free and $|T| \geq \frac{2h}{c}$ then it contains a smooth $(c,\lambda)$-$l$-sequence of length
$u$, where $c=(\frac{\lambda^{h}}{4^{h+3}h^{2}u^{h}(\lceil\log(u) \rceil)^{2h}})^{\lceil \log(u) \rceil}$.
\end{lemma}

\Proof
Take an arbitrary $0 < \lambda_{2} < 1$.
Without loss of generality we will assume that $u=2^{b}$ for some integer $b \geq 0$.
Take $n=|T|$. 
If $\lfloor \frac{cn}{h} \rfloor \geq \frac{cn}{2h}$ then
by the previous lemma we know that $V(T)$ contains two disjoint sets $X,Y$, each of size:
$c_{2}n$, where: $c_{2} = \frac{\lambda_{2}^{h}c}{h}$, where $c=\frac{1}{2h}$ and such that 
$d(X,Y) \geq 1 - \lambda_{2}$. Assume now that in every $H$-free tournament $T$ of order at least 
$n_{b-1}$ and for every $0 < \lambda_{b-1} < 1$ one can find in $T$ an $(c_{b-1},\lambda_{b-1})$-$l$-sequence of length $u^{'}=2^{b-1}$.
If this is the case then we can find in $X$ one $(\frac{c_{2}}{2}c_{b-1},\lambda_{b-1})$-$l$-sequnce
$Seq^{1}_{1}$ and remove it, then the next $(\frac{c_{2}}{2}c_{b-1},\lambda_{b-1})$-$l$-sequnce $Seq^{1}_{2}$
and remove it and so on...We can continue the procedure as long as we have at least $\frac{|X|}{2}$
vertices left. An analogous procedure can be applied to $Y$ to obtain 
$(\frac{c_{2}}{2}c_{b-1},\lambda_{b-1})$-$l$-sequnces $Seq^{2}_{1}, Seq^{2}_{2}...$.

We also need to assume that the size of the tournament from which an $l$-sequence is excluded is
big enough. We will get back to this assumption later while deriving the lower bound on the order
of $T$ from the assumptions of the lemma.
Denote by $X_{1}$ a subset of $X$ created by combining all $l$-sequences $Seq^{1}_{1},Seq^{1}_{2},...$
and by $Y_{1}$ a subset of $Y$ created by combining all $l$-sequences $Seq^{2}_{1},Seq^{2}_{2},...$
We have $|X_{1}| \geq \frac{|X|}{2}$, $|Y_{1}| \geq \frac{|Y|}{2}$. Therefore a simple
density argument gives us: $d(X_{1},Y_{1}) \geq 1 - 4\lambda_{2}$.
That means in particular that there exists $i,j$ such that 
$d(V(Seq^{1}_{i}),V(Seq^{2}_{j})) \geq 1 - 4\lambda_{2}$, where:
$V(Seq^{1}_{i}), V(Seq^{2}_{j})$ stand for the sets of vertices of the $l$-sequences
$Seq^{1}_{i}$ and $Seq^{2}_{j}$.
Now, from what we have said so far, we easily see that if we combine these 
two $l$-sequences $Seq^{1}_{i}$ and $Seq^{2}_{j}$ we get an $(c_{b},\lambda_{b})$-$l$-sequence,
where: $c_{b}=\frac{c_{2}}{2}c_{b}$ and $\lambda_{b}=\max(4\lambda_{2}(b-1)^{2},\lambda_{b-1})$
(again, by a simple density argument). As in one of the previous lemmas, we can easily
extract from it a smooth $l$-sequence that is a $(\frac{c_{b}}{2},4u\lambda_{b})$-$l$-sequence
(we leave details to the reader since the analysis is completely analogous to the one presented
earlier).
Thus we have: $c_{2}=\frac{\lambda_{2}^{h}}{2h^{2}}$ and $c_{b}=\frac{c_{2}}{2}c_{b-1}$,
$\lambda_{b}=\max(4\lambda_{2}(b-1)^{2},\lambda_{b-1})$ for $b>2$.
Solving this recurrence we get: $c_{b}=(\frac{c_{2}}{2})^{b}$, $\lambda_{b}=4\lambda_{2}(b-1)^{2}$.
Thus we get: $c_{b}=(\frac{\lambda_{2}^{h}}{4h^{2}})^{\log(u)}$.
It suffices to have: $4u \cdot 4\lambda_{2}\log(u)^{2} = \lambda$, i.e.: 
$\lambda_{2}=\frac{\lambda}{16u\log(u)^{2}}$. Substituting this expression on $\lambda_{2}$ into 
the formula for $c_{b}$ gives us: $c_{b}=(\frac{\lambda^{h}}{4^{h+2}h^{2}u^{h}(\log(u))^{2h}})^{\log(u)}$.
It remains to notice that our analysis is valid if $\lfloor \frac{c_{b}n}{h} \rfloor \geq
\frac{c_{b}n}{2h}$. Thus it suffices to have: $n \geq \frac{2h}{c_{b}}$. That completes the proof.

\bbox

We are ready to finish the proof of Theorem~\ref{conjecturetheorem}.
Notice that we procced by induction on $|T|$.
Take $\lambda=\frac{1}{1536 t^{2}(t+1)h^{3}2^{(t+1)h^{2}+h}}$.
We can now use Lemma~\ref{finallemma} and conclude that T contains a smooth
$(c,\lambda)$-$l$-sequence of length $2t+1$ for $c = \Omega(e^{-h^{4}\log(h)})$ for
$|T| \geq \frac{2h}{c}$. Now we can use Lemma~\ref{lintotrans} and extract from this
$l$-sequence an $(\frac{1}{4},\frac{c}{2},\frac{1}{24h^{3}2^{(t+1)h^{2}+h}},\epsilon)$
for $\epsilon = \Omega(\frac{1}{h^{4}\log(h)})$. Now we use Lemma~\ref{msequencelemma}.

We have: $M_{2} \leq 16h$ and $M_{1} = O(e^{h^{4}\log(h)})$. Thus we get: 
$\epsilon = \frac{\log(1 + \frac{1}{M_{2}-1})}{\log(M_{1})} = \Omega(\frac{1}{(M_{2}-1)\log(M_{1})})$.
Thus we get: $\epsilon = \Omega(\frac{1}{h^{5}\log(h)})$. In order to use Lemma~\ref{msequencelemma}
we need to have: $|T| \geq \frac{2h}{\frac{1}{4}}$ and $(1-\frac{2eh}{cn})^{h} \geq \frac{1}{2}$.
It is easy to see that all the lower bounds on $T$ we need to use all three lemmas are trivially 
satisfied for $|T| \geq W e^{h^{5}}$ for sufficiently large constant $W>0$.
On the other hand, for every fixed $W>0$ the theorem is trivially true for 
$\epsilon=\frac{w}{h^{5}\log(h)}$ for small enough constant $w>0$
and all tournaments with at most $W e^{h^{5}}$ vertices. That completes the proof.

\bbox

Now we will prove Theorem~\ref{algorithmtheorem}.

\Proof
The coloring algorithm extracts big transitive subtournaments one be one as long there are some vertices left in the $H$-free tournament $T$.
To extract big transitive subtournaments it uses Theorem~\ref{conjecturetheorem}. It is easy to see that such a procedure produces a partitioning of the
vertices of $T$ into at most $n^{1-\epsilon}\log(n)$ transitive subtournaments, where  $\epsilon$ is the lower bound on $\epsilon(H)$ as in Theorem~\ref{conjecturetheorem}.
To upper-bound the running time of this approach we need tol find an upper-bound on the running time of the subroutine finding big transitive subtournament in the $H$-free tournament $T$.
We will use the proof of the Theorem~\ref{conjecturetheorem}. Note that almost all the steps of the proof of  Theorem~\ref{conjecturetheorem} can be directly
translated into polynomial subroutines. There are two exceptions: the part where transitive subtournaments are being inductively extracted from linear sets and the
part where $l$-sequences are being inductively extracted from linear sets. That observation and the analysis of these two parts easily  lead to the following recursive
 formula on the total running time $T(n)$ of the algorithm (we leave details to the reader): $T(n) \leq T((1-c)n) n^{1-\epsilon} + poly(n)$, where $c$ is some constant (parameter not depending on the size of the $H$-free tournament $T$) .
Solving this recurrence gives us the upper bound on the running time as in Theorem~\ref{algorithmtheorem}.
\bbox

\section{Proof of Theorem~\ref{generaltheorem} }

We are now ready to prove Theorem~\ref{generaltheorem}.

\Proof
To prove that $\epsilon(H)  = \Omega({e^{-(|H|+1)\log(|H|+1)}})$ for every known EH tournament 
we use the substitution procedure and related theoretical guarantees for the EH coefficient of the 
outcome tournament (see: \cite{alon}) as well as our earlier result.
Let $f$ be a nondecreasing function taking positive values.
Assume that one can prove that $\epsilon(H) \geq \frac{1}{f(|H|)}$ for all known EH tournaments $H$ 
of size at most $r$ and all known EH prime tournaments. Let $H$ be a known EH prime tournament of
size $r+1$. Assume that it is not prime. Thus it can be constructed from some tournaments: $D,F$ of smaller
orders by replacing one vertex of $D$ with the copy of $F$, according to the substitution procedure.
Following \cite{alon}, we obtain: 

\begin{equation}
\label{alonformula}
\epsilon(H) \geq \frac{\epsilon(F)\epsilon(D)}{\epsilon(D)+k \epsilon(F)} - \epsilon
\end{equation}
for every $\epsilon > 0$.
Denote: $|D|=k, |F|=l$. Then we have: $|H| = k + l -1$.
To prove that $\epsilon(H) \geq \frac{1}{f(|H|)}$ it suffices (according to inequality \ref{alonformula}) to prove that:
$\frac{\frac{1}{f(l)}\frac{1}{f(k)}}{\frac{1}{f(k)} + \frac{k}{f(l)}} > \frac{1}{f(k+l-1)}$, i.e. that:
$f(k+l-1) > f(l) + kf(k)$. Denote: $f(i) = e^{t(i)}$, where $t$ is a nondecreasing function.
We want: $e^{t(k+l-1)} \geq e^{t(l)} + e^{\log(k) + t(k)}$.
Notice first that we have: $k,l \geq 2$. Assume first that $k \geq l$. Then it is easy to see that
it suffices to have: $t(k+l-1) \geq t(k) + \log(k) + \log(2)$. Under assumption that $t$ is nondecreasing
we see that it suffices to have: $t(k+1) \geq t(k) + \log(k) + \log(2)$. Notice that trivially we can take:
$t(1)=0$. 

Thus we conclude that it is enough to take: $t(k) = \log(k-1)! + (k-1)\log(2)$.
So it suffices to have: $t(k) = k \log(2k)$. Let us assume now that $l \geq k$.
The we want to get: $e^{t(k+l-1)} \geq e^{t(l)}(1 + e^{t(k)-t(l)+\log(k)})$.
It is enough to have: $e^{t(l+1)} \geq e^{t(l)}(1 + e^{t(k)-t(l)+\log(k)})$.
Taking $t$ to be nondecreasing it suffices to get: $e^{t(l+1)} \geq e^{t(l)}(1+l)$.
Thus it is enough to have: $e^{t(l+1)} = (l+1)!$, i.e. $t(l+1) = \log(l+1)!$.
Theorem \ref{conjecturetheorem} says that there exists a constant $C>0$ such that every prime tournament 
$H$ satisifes: $\epsilon(H) >= \frac{C}{|H|^{5}\log(|H|)}$. Combining this with the analysis of function $t$,
we conclude that every tournament $H$ that can be obtained from known prime EH tournaments by 
the substitution procedure satisfies: $\epsilon(H) \geq \frac{C}{e^{(|H|+1)\log(|H|+1)}}$.
Since all known EH tournaments are those that are constructed from known prime EH tournaments by the 
substitution procedure, we are done.
\bbox

\section{Proof of Theorem~\ref{uppertheorem}}

In this section we give the proof of Theorem~\ref{uppertheorem}.
First we summarize previously existing results regarding upper bounds on EH coefficients of tournaments.
The first upper bounds on the EH coefficients of random tournaments were given in \cite{kchoromanski},
where it was proven that:

\begin{theorem}
 There exists $\eta>0$ such that if $\mathcal{H}^{n,\eta}$ denotes the set of $n$-vertex tournaments
 satisfying $\epsilon(H) \leq \frac{4}{n}(1+\frac{\sqrt{\log(n)}}{\sqrt{n}})$ and $\mathcal{H}^{n}$
 denotes the set of all $n$-vertex tournaments then:
 $$\lim_{n \to \infty} \frac{|\mathcal{H}^{n,\eta}|}{|\mathcal{H}^{n}|} = 1.$$
\end{theorem}

In other words, random tournaments have EH coefficients of the order $O(\frac{1}{|H|})$.

Surprisingly, it turns out that the partitioning number $p(H)$ of a tournament $H$ tells 
us something about EH coefficient of $H$.
It is a well-known fact (see: Appendix A) that:

\begin{theorem}
 There exists $C>0$ such that $\epsilon(H) <= C\frac{\log(\log(p(H)))}{\log(p(H))}$.
\end{theorem}

Besides, under assumption that $H$ is prime, that result was strengthened.  It is known (see: Appednix B)
that:

\begin{theorem}
 There exists $C>0$ such that every prime tournament $H$satisfies:
 $\epsilon(H) \leq C\frac{\log(|H|)}{|H|}$.
\end{theorem}

Now lets look on the upper bounds for EH coefficients written as functions of tournaments' partitioning
numbers. Note that, since for prime tournaments $H$ we have: $p(H)=|H|$, we get the following
bounds on EH coefficients for prime tournaments $H$: $\epsilon(H) \leq C\frac{\log(p(H))}{p(H)}$. 
There is a striking difference between an expression for the upper bound on the EH coefficient 
for a prime tournament $H$ and for a general tournament that does not have to be necessarily prime.
It seems that bounds given for general tournaments $H$ can be significantly improven.
This is in fact true, at least if a tournament does not have too large homogeneous sets.
Our improvement led to polynomial lower and upper bounds on EH coefficients for several classes
of known EH tournaments. Those much tighter bounds is a step towards understanding how the EH
coefficients depend on the order of $H$.

We are ready to prove Theorem~\ref{uppertheorem}.

\Proof

Denote $h=|H|$. We can assume that $h>1$. Let $B$ be a tournament and denote: $n=|B|$. 
Following the procedure described in Appendix A,
we define the family of tournaments: $T^{B}_{0},T^{B}_{1},...$ as follows:

\begin{itemize}
 \item $T^{B}_{0}$ is a single vertex,
 \item $T^{B}_{k+1}$ is obtained from $T^{B}_{k}$  by replacing each vertex $v \in V(B)$ with
       the copy $T^{v}$ of $T^{B}_{k}$ and making a vertex $u_{1} \in T^{v}$ adjacent to a vertex
       $u_{2} \in T^{w}$ iff $v$ is adjacent to $w$ in $B$ for $k=0,1,2,...$.
\end{itemize}

It was shown in Appendix A that $T^{B}_{k}$ does not contain transitive subtournaments of
order larger than $|T^{B}_{k}|^{\epsilon}$, where: $\epsilon = \frac{\log(t)}{\log(n)}$ and
$t$ is the largest transitive subtournament of $B$.
Denote by $\mathcal{Q}^{H}$ the family of all quotients tournaments of a tournament $H$
of order greater than one and by $\mathcal{Q}^{H}_{i}$ the family of all quotient tournaments
of a tournament $H$ of exactly $i$ vertices.
Note that from the definition of the partitioning number it follows that:
\begin{equation}
 \mathcal{Q}^{H} = \mathcal{Q}^{H}_{p(H)} \cup ... \cup \mathcal{Q}^{H}_{|H|}.
\end{equation}
Note now that if $B$ does not contain any tournament from 
$\mathcal{Q}^{H}_{p(H)} \cup ... \cup \mathcal{Q}^{H}_{|H|}$ as a subtournament then
every $T^{B}_{k}$ is $H$-free. Therefore to finish the proof of the theorem it suffices to
construct for every $\eta>0$ a tournament $B$ that is $W$-free for every 
$W \in \mathcal{Q}^{H}_{p(H)} \cup ... \cup \mathcal{Q}^{H}_{|H|}$
and does not contain transitive subtournaments of order larger than
$n^{\rho}$, where $\rho=C(\eta) \frac{\log(p(H))}{p(H)^{\frac{1}{2}-\eta}}$ and
$C(\eta)$ is independent of $H$.
Denote $h=|H|$. Since the biggest homogeneous set of $H$ has no more than 
$\frac{\sqrt{h}}{2}$ vertices we conclude that 
$p(H) \geq \frac{h}{\frac{\sqrt{h}}{2}}$. Thus we have: $p(H) \geq \sqrt{h}$.
Our goal is to find an upper bound $f(i)$ on $|\mathcal{Q}^{H}_{i}|$ that does not depend
on $H$ but only on $i$. Before doing it we solve a little bit simpler task.
Denote by $\mathcal{Q}^{\epsilon_{2}}_{i}$ the family of $i$-vertex 
tournaments that contain homogeneous set $\mathcal{S}$ satisfying 
$i^{\epsilon_{2}} \leq |\mathcal{S}| \leq \frac{i}{2}$, 
for $\epsilon_{2}=\frac{1}{2} - \eta$ (note that without loss of generality we an assume that $\eta < \frac{1}{2}$).
To find the upper bound on 
$|\mathcal{Q}^{\epsilon_{2}}_{i}|$ we use the probabilistic argument.
Let $R$ be a random tournament on $i$ vertices where the direction of each edge
is chosen independently at random and each direction has probability $\frac{1}{2}$ of
being chosen. Denote by $X$ the random variable that counts the number of
homogeneous sets in $R$ of size at least $i^{\epsilon_{2}}$ and at most $\frac{i}{2}$.
For every fixed subset $\mathcal{S} \subseteq V(R)$ satisfying 
$i^{\epsilon_{2}} \leq |\mathcal{S}| \leq \frac{i}{2}$ we know that the probability $p$ that it is
homogeneous is at most $(\frac{2}{2^{i^{\epsilon_{2}}}})^{\frac{i}{2}}$.
This comes directly from the fact that each point from $V(R) \backslash \mathcal{S}$ is either
adjacent to all the vertices of $\mathcal{S}$ or adjacent from all the vertices of $\mathcal{S}$,
the size of $V(R) \backslash \mathcal{S}$ is at least $\frac{i}{2}$ and $|\mathcal{S}| \geq i^{\epsilon_{2}}$.
Thus we have:
\begin{equation}
 EX \leq 2^{i} (\frac{2}{2^{i^{\epsilon_{2}}}})^{\frac{i}{2}}. 
\end{equation}
Thus $P(X>0) \leq EX \leq 2^{\frac{3i}{2} - \frac{1}{2}i^{1+\epsilon_{2}}}$.
We can conclude that $|\mathcal{Q}^{\epsilon_{2}}_{i}| \leq 2^{{i \choose 2}} 2^{\frac{3i}{2} - \frac{1}{2}i^{1+\epsilon_{2}}}$.
We have now two possibilities:
\begin{itemize}
\label{poss}
 \item $|\mathcal{Q}^{H}_{i}| \leq 2^{{i \choose 2}} 2^{\frac{3i}{2} - \frac{1}{2}i^{1+\epsilon_{2}}}$, or
 \item $|\mathcal{Q}^{H}_{i}| > 2^{{i \choose 2}} 2^{\frac{3i}{2} - \frac{1}{2}i^{1+\epsilon_{2}}}$.
\end{itemize}

Assume first that the latter holds.
In the latter case 
we have $\mathcal{Q}^{H}_{i} \backslash \mathcal{Q}^{\epsilon_{2}}_{i} \neq \emptyset$
thus $\mathcal{Q}^{H}_{i}$ contains a tournament  such that its homogeneous sets are of size smaller than 
$i^{\epsilon_{2}}$ or larger than $\frac{i}{2}$.
If it contains a homogeneous set of size larger than $\frac{i}{2}$ then $H$ contains a homogeneous set of size
at least $\frac{i}{2}$. Note that $i \geq p(H)$ and since we have already
proved that $p(H) \geq \sqrt{h}$, we concldue that $H$ contains a homogeneous
set of size at least $\frac{\sqrt{h}}{2}$ which is a contradiction.
Thus we have in $\mathcal{Q}^{H}_{i}$ a tournament $W$ with the largest
homogeneous set of size smaller than $i^{\epsilon_{2}}$.
Denote $V(W)=\{w_{1},...,w_{i}\}$.
Let $L$ be a tournament from $\mathcal{Q}^{H}_{i}$ and denote
$V(L)=\{l_{1},...,l_{i}\}$. Notice that each vertex of $V(W)$ is related to the
homogeneous set from some homogeneous partitioning of $V(W)$.
The same is true about the vertices of $V(L)$.
Denote the homogeneous partition related to $W$ as:
$\{W_{1},...,W_{i}\}$ and the homogeneous partition related to $L$ as
$\{L_{1},...,L_{i}\}$. We assume that $w_{j}$ corresponds to $W_{j}$
and $l_{j}$ corresponds to $L_{j}$ for $j=1,2,...,i$.
We will now define the bipartite graph $G_{W,L}$ as follows:
\begin{itemize}
\item $V(G_{W,L}) = V(W) \cup V(L)$,
\item the color classes of $G(W,L)$ are: $V(W)$ and $V(L)$,
\item there is an edge between $w_{j_{1}}$ and $l_{j_{2}}$ if and only if
      $W_{j_{1}} \cap L_{j_{2}} \neq \emptyset$.
\end{itemize}

\begin{lemma}
\label{forbiddenlemma}
 $G_{W,L}$ does not contain a cycle of length four $C_{4}$ as a subgraph.
\end{lemma}

\Proof
Assume by contradiction that there exist 
$w_{k_{1}},w_{k_{2}}, l_{s_{1}},l_{s_{2}}$ for some $k_{1} < k_{2}$,
$s_{1} < s_{2}$, that induce $C_{4}$ in $G_{W,L}$.
Take corresponding subsets $W_{k_{1}},W_{k_{2}},L_{s_{1}},L_{s_{2}}$.
Let $a \in W_{k_{1}} \cap L_{s_{1}}$, $b \in W_{k_{1}} \cap L_{s_{2}}$,
$c \in W_{k_{2}} \cap L_{s_{1}}$, $d \in W_{k_{2}} \cap L_{s_{2}}$.
Assume without loss of generality that $W_{k_{1}}$ is adjacent to $W_{k_{2}}$.
Thus we have: $a$ is adjacent to $d$. Similarly, $b$ is adjacent to $c$.
But then we get a contradiction since $\{a,c\} \subseteq L_{s_{1}}$,
$\{b,d\} \subseteq L_{s_{2}}$ and either $L_{s_{1}}$ is adjacent to $L_{s_{2}}$
or vice versa.

\bbox

A classic result in extremal graph theory states that the number of bipartite
graphs with color classes of size $i$ is at most $2^{c_{1} i^{\frac{3}{2}}}$
for some universal constant $c_{1}>0$.

Fix graph $W \in \mathcal{Q}^{H}_{i}$ with the largest homogeneous set of size smaller 
than $i^{\epsilon_{2}}$. Using our previous remarks, we conclude that the number of 
different bipartite graphs $G_{W,L}$ is at most $2^{c_{1} i^{\frac{3}{2}}}$.
Fix bipartite graph $G_{W,L}$. Take two vertices $l_{s_{1}}, l_{s_{2}} \in V(L)$.
Note that type of adjacency between $l_{s_{1}}$ and $l_{s_{2}}$ is uniquely determined
by the graph $G_{W,L}$ unless both $L_{s_{1}}$ and $L_{s_{2}}$ are subsets of some $W_{k}$.
Let us count the maximal number of vertices of $L$ such that all corresponding
sets $L_{j}$ are subsets of some $W_{k}$.
Denote $d_{0}=\frac{i-1}{i^{\epsilon_{2}}}$. If there exists $W_{k}$ and different sets $L_{t_{1}},...,L_{t_{m}}$
such that $m \geq i - d_{0}$ then, by Pigeonhole Principle, there exists 
$L_{j}$ for some $j \in \{t_{m}+1,...L_{i}\}$ that intersects at least $\frac{i-1}{d_{0}}$
sets from $\{W_{1},...,W_{i}\}$.
Denote those sets as $W_{p_{1}},...,W_{p_{u}}$, where $u \geq \frac{i-1}{d_{0}}$.
Now take some $W_{z} \notin \{W_{p_{1}},...,W_{p_{u}}\}$. Set $W_{z}$ intersects with some
set from $\{L_{1},...,L_{i}\}$ but we already now that this is not $L_{j}$.
Denote this set by $L_{b}$. Then, if $L_{b}$ is adjacent to $L_{j}$ then $W_{z}$ is adjacent to 
every $W_{p_{i}}$. Similarly, if $L_{j}$ is adjacent to $L_{b}$ then every $W_{p_{i}}$ is
adjacent to $W_{z}$. Because $W_{z}$ was chosen arbitrarily, we conclude that 
$\{w_{p_{1}},...,w_{p_{u}}\}$ is a homogeneous set in $V(W)$.
Since $u > \frac{i-1}{d_{0}}$ and from our choice of $d_{0}$ we obtain that $W$ contains 
a homogeneous set of size more than $i^{\epsilon_{2}}$. This however contradicts the 
definition of $W$. We conclude that every $W_{i}$ contains at most $i-d_{0}$ sets from
$\{L_{1},...,L_{i}\}$ as subsets.

Let us conclude what we have managed to show so far.
We know that either 
$|\mathcal{Q}^{H}_{i}| \leq 2^{{i \choose 2}} 2^{\frac{3i}{2} - \frac{1}{2}i^{1+\epsilon_{2}}}$
or there exists $W \in \mathcal{Q}^{H}_{i}$ such that:
\begin{itemize}
 \item the number of all bipartite graphs $G_{W,L}$, where $L \in \mathcal{Q}^{H}_{i}$
       is at most $2^{c_{1}i^{\frac{3}{2}}}$,
 \item no $W_{i}$ from the homogeneous partitioning corresponding to $W$ contains more than
       $i-d_{0}$ sets $L_{j}$ from the partitioning corresponding to $L$,
 \item for $l_{s_{1}}, l_{s_{2}} \in V(L)$ type of adjacency between $l_{s_{1}}$ and 
           $l_{s_{2}}$ is determined by $G_{W,L}$ unless both $L_{s_{1}}$ and $L_{s_{2}}$
           are subsets of some $W_{j}$.
\end{itemize}
We conclude that as an upper bound $f$ on $|\mathcal{Q}^{H}_{i}|$ we can take any function
$N \rightarrow R$  
satisfying the following three properties (that we will call \textit{basic properties}):
\begin{itemize}
\item  $f(i)$ is a valid bound for $i$ small enough,
\item  $f(i) \geq 2^{c_{1}i^{\frac{3}{2}}} \max_{k_{1},...,k_{i}} f(k_{1}) \cdot .... \cdot f(k_{i})$ for larger $i$,
\item  $\forall_{i} f(i) \geq 2^{{i \choose 2}} 2^{\frac{3i}{2} - \frac{1}{2}i^{1+\epsilon_{2}}}$,
\end{itemize}

where $k_{1},....,k_{i} \in N$, $k_{1}+...+k_{i}=i$ and $k_{1},...,k_{i} \leq i-d_{0}$.

We will take $f>0$ such that $\log(f)$ is convex in the domain $[1, \infty]$
and $\log(f(1))-\log(f(0)) \leq \log(f(2))-\log(f(1))$.
We call such a function an \textit{$\alpha$-function}. 
Let assume now that $f$ is an $\alpha$-function.
Take $k_{1},k_{2} \geq 0$ such that $k_{2} \geq 1$ and $k_{1} \geq k_{2}$.
We have: 
\begin{equation}
 \log(f(k_{1}+1))-\log(f(k_{1})) \geq \log(f(k_{2})) - \log(f(k_{2}-1)).
\end{equation}

The inequality above comes directly from the definition of an $\alpha$-function.
Thus we have:
\begin{equation}
\label{impen}
 f(k_{1})f(k_{2}) \leq f(k_{1}+1)f(k_{2}-1).
\end{equation}

Let us take $f(i) = 2^{ci + {i \choose 2} - \frac{1}{2}\delta(i)}$, where $\delta(i) = i^{\frac{3}{2}-\eta}$
and $c$ is a big enough constant (without loss of generaliy we can assume that $\eta$
is smaller than any fixed in advance constant, in particular: $\eta < \frac{3}{2}$).
For $c$ big enough $f(i)$ is a valid bound for every $i \leq M$, where $M$ is any
constant chosen in advance. Thus first basic property is satisfied. The last one is clearly satisfied too.
Note also that it is easy to check that $f$ is an $\alpha$-function.
Now observe that, according to inequality~\ref{impen}, to prove that the second basic property is satisfies
it only suffices to
show that:
\begin{equation}
\label{ineq2}
 f(i) \geq 2^{c_{1}i^{\frac{3}{2}}}f(i-d_{0})f(d_{0})f(0)^{i-2}
\end{equation}

for $i > M$. Using our formula on $f$ and carefully evaluating one can easily check
that $f$ defined above satisfies inequality~\ref{ineq2} for $i>M$, where $M$ large ebough.
Thus we can conclude that there exists $c>0$ such that:
$$|\mathcal{Q}^{H}_{i}| \leq f(i),$$
where $f(i) = 2^{ci + {i \choose 2} - \frac{1}{2}\delta(i)}$ and $\delta(i) = i^{\frac{3}{2}-\eta}$.

Now let us construct base tournament $B$ with $n = e^{c_{3} p(H)^{\frac{1}{2}-\eta}}$, where
$c_{3}>0$ is small enough and such that 
directions of edges are chosen independently at random (each of two possibilities with probability
$\frac{1}{2}$). Let $Z$ be a random variable that counts the number of copies of tournaments
from $\mathcal{Q}^{H}$ that are in $B$.
Note that we have:
$$EZ \leq \sum_{i=p(H)}^{|H|} i! {n \choose i} \frac{1}{2^{{i \choose 2}}} f(i).$$
Evaluating that expression, one can easily check that 
$EZ \leq \sum_{i=p(H)}^{|H|} e^{-c_{4}i^{\frac{3}{2}-\eta}}$, where $c_{4}>0$ is some constant.
Since the sum $\sum_{i=1}^{\infty} e^{-c_{4}i^{\frac{3}{2}-\eta}}$ is finite, we conclude that
$EZ \leq C$ for some constant $C>0$. Thus by deleting from $B$ at most $C$ vertices we
can make it $Q$-free for every $Q \in \mathcal{Q}^{H}$.
Therefore to complete the entire proof it suffices to prove that: 
$\frac{\log(t)}{\log(n-C)} \leq C(\eta)\frac{\log(p(T))}{p(T)^{\frac{1}{2}-\eta}}$, where
$t$ is the size of the largest transitive subtournament of $B$. But since $B$ is random,
with overwhelming probability all its transitive subtournaments are of order at most $\log(n)$
and so the inequality follows.
\bbox

\bibliographystyle{unsrt}
\small{
\bibliography{galaxy_coloring}
}

\appendix

\section{State-of-the-art techniques regarding upper bounds on EH coefficients}

Before proving Theorem \ref{uppertheorem} we will present some well-known results regarding deriving upper
bounds on EH coefficients of tournaments that we refered/will be refering to. 
These results were sent for publication but were not published yet.
Therefore, for completeness we will give here
their full statements and proofs. We want to emphasize that these tools do not stand for the contribution we make in the paper.

We will use two results:

\begin{theorem}
\label{homotheorem}
There exists $C>0$ such that every prime $h$-vertex tournament $H$ satisfies $\xi(H) \leq C \frac {\log(h)}{h}$. 
\end{theorem}

\begin{theorem}
\label{stronghomotheorem}
There exists $C>0$ such that every $h$-vertex tournament $H$ satisfies $\xi(H) \leq C \frac{\log(\log(p(H)))}{\log(p(H))}.$ 
\end{theorem}

We need one more definition.

\begin{definition}
 For a tournament $H$ we say that $Q$ with $V(Q)=\{q_{1},...,q_{|Q|}\}$ is a 
 \textit{quotient tournament of $H$} if there exists
 a homogeneous partitioning $\{V_{1},...,V_{|Q|}\}$ of the vertices of $V(H)$ such that
 $V_{i}$ is adjacent to $V_{j}$ iff $q_{i}$ is adjacent to $q_{j}$ in $Q$ 
 (thus $V_{j}$ is adjacent to $V_{i}$ iff $q_{j}$ is adjacent to $q_{i}$ in $Q$).
\end{definition}

\subsection{Proofs of Theorem \ref{homotheorem} and Theorem \ref{stronghomotheorem}}

Denote by $tr(T)$ the largest size of the transitive subtournament of a tournament $T$. For $X \subseteq V(T)$, write $tr(X)$ for $tr(T|X)$. Let $H$ be a tournament.
Assume that $V(H)$ admits a homogeneous partitioning $P=\{V_{1},...,V_{k}\}$. We associate with the partitioning $P$ a $k$-vertex quotient tournament $H_{P}$ with $V(H_{P})=\{v_{1},v_{2},...,v_{k}\}$ such that for $1 \leq i < j \leq k$ vertex $v_{i}$ is adjacent to a vertex $v_{j}$ in $H_{P}$ if $V_{i}$ is complete to $V_{j}$. 
We say that a tournament $T$ is $H$-far if $T$ is $H_{P}$-free for every $k>1$ and every homogeneous partitioning $P$ of $H$ consisting of $k$ parts.

First we prove the following result:
\begin{theorem}
\label{strongerhomotheorem}
Let $H$ be a tournament with at least two vertices. Assume that $T$ is $H$-far. Then $$\xi(H) \leq \frac{\log(tr(T))}{\log(|T|)}.$$
\end{theorem}

\Proof
Denote $V(T)=\{1,2,...,|T|\}$.
Consider a family of tournaments $\{F_{0},F_{1},...\}$ defined in the following recursive way.
A tournament $F_{0}$ is just a single vertex. For $i > 0$ a tournament $F_{i}$ is defined as follows.
$V(F_{i})=P^{i}_{1} \bigcup P^{i}_{2} \bigcup ... P^{i}_{|T|}$, where each $P^{i}_{j}$ for $j=1,2,...,|T|$ induces a tournament isomorphic to $F_{i-1}$ and besides for any two $1 \leq j_{1} < j_{2} \leq |T|$ the set $P^{i}_{j_{1}}$ is complete to the set $P^{i}_{j_{2}}$ if $j_{1}$ is adjacent to $j_{2}$ in $T$ and complete from the set $P^{i}_{j_{2}}$ if $j_{1}$ is adjacent from $j_{2}$.  
Note first that every $F_{i}$ is $H$-free. To see this we use induction on $i$. For $i=0$ this is trivial.
Now take tournament $F_{i+1}$. If $F_{i+1}$ is not $H$-free, then since $T$ is $H$-far and $F_{i}$ is $H$-free, we can conclude that $V(H)$ has homogeneous partitioning consisting of $k$ parts for some $1< k < p(H)$. That contradicts definition of $p(H)$. 
Knowing that every $F_{i}$ is $H$-free we calculate the size of the biggest transitive subtournament of $F_{i}$. For $i=0$ we have $tr(F_{i})=1$. Assume that $i>0$. Let $Tr_{i}$ be the biggest transitive subtournament of $F_{i}$. Write $S_{j} = V(Tr_{i}) \bigcap P^{i}_{j}$ for $j=1,2,...,|T|$. Assume that $\{S_{j_{1}},...,S_{j_{k}}\}$ is the set of nonempty sets $S_{j}$. Note that the subtournament of $T$ induced by the set $\{j_{1},...,j_{k}\}$ must be transitive. Otherwise, according to the definition of the family $\{F_{j}\}_{j=0,1,2,...}$, we  conclude that $Tr_{i}$ contains vertices inducing directed triangle (that contradicts the fact that $Tr_{i}$ is transitive).
Therefore we must have $k \leq tr(T)$. Since $S_{j} \subseteq P^{i}_{j}$ we must have $|S_{j}| \leq tr(F_{i-1})$. Therefore we have $|V(Tr_{i})|=tr(F_{i}) \leq tr(T)tr(F_{i-1})$.
So by induction, $tr(F_{i}) \leq tr(T)^{i}$. In fact from our analysis we easily see that we have
$tr(F_{i})=tr(T)^{i}$. We also have $|V(F_{i})|=|T|^{i}$. Therefore we have $tr(F_{i}) = |F_{i}|^{\log_{|T|}(tr(T))}$. So we have $tr(F_{i}) = |F_{i}|^{\frac{\log(tr(T))}{\log(|T|)}}$. We conclude that each $F_{i}$ is $H$-free and does not contain transitive subtournaments of size at least $|F_{i}|^{\epsilon}$, where $\epsilon=\frac{\log(tr(T))}{\log(|T|)}$. This implies that $\xi(H) \leq \epsilon$.
\bbox

We are now ready to prove Theorem \ref{homotheorem} and Theorem \ref{stronghomotheorem} that we encapsulate in the following statement:

\begin{theorem}
\label{unitheorem}
There exists $C>0$ such that every $h$-vertex tournament $H$ satisfies $\xi(H) \leq C \frac{\log(\log(p(H)))}{\log(p(H))}.$ Furthermore, if $p(H)=h$ then $\xi(H) \leq C \frac{\log(h)}{h}$ for some universal constant $C>0$.
\end{theorem}

\Proof 

We may assume that $p(H)$ is large enough since for every tournament $H$ we trivially have: $\xi(H) \leq 1$.
Let $G$ be a $n$-vertex tournament, where for any two vertices $1 \leq i < j \leq n$ an edge $(i,j)$ is chosen with probability $\frac{1}{2}$. Let $c$ be some large constant. Denote by $X$ the number of  transitive subtournaments of $G$ of size at least $c\log(n)$ and by $Y$ the number of copies in $G$ of subtournaments isomorphic to some $H_{P}$, where $P$ is some homogeneous partitioning of $H$. Write $r=c\log(n)$. Note that we have $EX \leq r!{n \choose r} (\frac{1}{2})^{{r \choose 2}} $.
Therefore $EX \leq e^{r\log(n)-\frac{r(r-1)}{2}\log(2)}$. Taking $c$ large enough we have $EX < \frac{1}{3}$.
Assume first that $p(H)=h$. Note that in this case there is a unique $H_{P}$ and it is isomorphic to $H$. Write $n=e^{dh}$, where $d>0$ is a small enough constant. We have:
$EY \leq {n \choose h} h! 2^{-{h \choose 2}} \leq n^{h} 2^{-{h \choose 2}} < \frac{1}{3}$ for $d$ small enough.
Therefore for $c$ large enough and $d$ small enough we have: $EX < \frac{1}{3}$ and $EY <\frac{1}{3}$. Thus, using Markov's inequality, we conclude that with probability less than $\frac{1}{3}$ we have $Y \geq 1$ and with probability less than $\frac{1}{3}$ we have $X \geq 1$.
So from the union bound we know that with probability bigger than $\frac{1}{3}$ we have $X < 1$ and $Y < 1$. So there exists a tournament $G$ that is $H$-far and does not contain transitvie subtournaments of size $c\log(n)$. Since we have $n=e^{p(H)}$, using Theorem~\ref{strongerhomotheorem}, we immediately obtain Theorem~\ref{homotheorem}.   
In the general case when the condition $p(H)=h$ is not necessarily satisfied, we use the same analysis. The only difference is the choice of $n$. Let $n=p(H)-1$. In this scenario $Y$ is trivially $0$ since every $H_{P}$ has at least $p(H)$ vertices so cannot be contained in the tournament of $p(H)-1$ vertices. The rest of the proof is exactly the same as in the case when $p(H)=h$.

\bbox

\end{document}